\documentclass[a4paper,reqno,twocolumn]{article}

\usepackage[utf8]{inputenc} 
\usepackage[T1]{fontenc} 
\usepackage{lmodern} 
\usepackage[american]{babel}
\usepackage{csquotes}

\usepackage{verbatim} 
\usepackage{authblk} 
\usepackage{lipsum}
\usepackage{calc}
\usepackage{amsmath}
\usepackage{amsfonts}
\usepackage{amssymb}
\usepackage{graphicx}
\usepackage{caption}
\usepackage{subcaption}
\usepackage{float}
\usepackage{xcolor}
\usepackage{enumitem}
\setlist[enumerate]{leftmargin=.5in} 
\setlist[itemize]{leftmargin=.5in} 
\providecommand{\keywords}[1]{\textbf{Keywords: } #1} 
\makeatletter 
\newcommand{\subjclass}[2][2020]{%
	\let\@oldtitle\@title%
	\gdef\@title{\@oldtitle\footnotetext{#1 \emph{Mathematics Subject Classification.} #2}}%
}
\makeatother
\usepackage{fancyhdr}
\allowdisplaybreaks 
\usepackage{scalerel,stackengine}
\stackMath
\newcommand\reallywidehat[1]{%
	\savestack{\tmpbox}{\stretchto{%
			\scaleto{%
				\scalerel*[\widthof{\ensuremath{#1}}]{\kern-.6pt\bigwedge\kern-.6pt}%
				{\rule[-\textheight/2]{1ex}{\textheight}}
			}{\textheight}%
		}{0.5ex}}%
	\stackon[1pt]{#1}{\tmpbox}%
}

\usepackage[linesnumbered,ruled,vlined]{algorithm2e}

\SetCommentSty{mycommfont}
\SetKwInput{KwInput}{Input}                
\SetKwInput{KwOutput}{Output}              

\usepackage[
pdfmenubar=true,
pdftitle={Strong fast invertibility and Lyapunov exponents},
pdfauthor={Florian Noethen},
colorlinks=true,
linkcolor=black,
hidelinks,
raiselinks=false,
hyperfootnotes=true,
linktoc=page
]{hyperref}
\usepackage[capitalize,nameinlink]{cleveref} 
\crefformat{footnote}{#2\footnotemark[#1]#3} 
\crefname{subsection}{Subsection}{Subsections}

\let\originalleft\left
\renewcommand{\left}{\mathopen{}\mathclose\bgroup\originalleft} 
\let\originalright\right\renewcommand{\right}{\aftergroup\egroup\originalright} 

\usepackage{amsthm}
\usepackage{amssymb}
\newtheoremstyle{break} 
{\topsep}{\topsep}%
{\itshape}{}%
{\bfseries}{}%
{\newline}{}%
\theoremstyle{break}
\newtheorem*{theorem*}{Theorem} 
\newtheorem{theorem}{Theorem}[section] 
\newtheorem{lemma}[theorem]{Lemma} 
\newtheorem{proposition}[theorem]{Proposition}
\newtheorem{definition}[theorem]{Definition}

\newtheorem{remark}[theorem]{Remark}
\newtheorem{corollary}[theorem]{Corollary}
\newtheorem*{conjecture*}{Conjecture} 

\Crefname{conjecture}{Conjecture}{Conjectures}

\usepackage[
	style = numeric,
	sorting = nyt,
	sortcites,
	giveninits = true,
	isbn = false,
	doi = true,
    url = true,
	clearlang = true,
	backend = biber,
	citestyle = numeric-comp  
]{biblatex}
\DeclareNameAlias{author}{family-given}
\AtEveryBibitem{\clearfield{pagetotal}} 
\AtEveryBibitem{\clearfield{note}} 
\AtEveryBibitem{\clearfield{month}}
\AtEveryBibitem{\clearfield{day}}
\DeclareSourcemap{
	\maps[datatype=bibtex]{
		\map{
			\pernottype{misc} 
			\step[fieldset=url, null] 
			\step[fieldset=urldate, null] 
		}
	}
}
\renewbibmacro{in:}{ 
	\ifentrytype{article}{}{\printtext{\bibstring{in}\intitlepunct}}%
}

\begin{document}
	
	\title{Strong Fast Invertibility and Lyapunov Exponents for Linear Systems}
\author{Florian Noethen\thanks{Universität Hamburg, Bundesstraße 55, 20146 Hamburg, Germany (\href{mailto:florian.noethen@uni-hamburg.de}{florian.noethen@uni-hamburg.de}, ORCID \href{https://orcid.org/0000-0001-9609-8206}{0000-0001-9609-8206}).}
		
}
\date{\today} 
\maketitle
\thispagestyle{empty}

\begin{abstract}
	
	In 2019 Anthony Quas, Philippe Thieullen and Mohamed Zarrabi introduced the concept of strong fast invertibility for linear cocycles. It relates the growth of volumes between different initial times and, together with a condition on singular value gaps, yields the existence of a dominated splitting of the dynamics. The properties of this splitting largely coincide with those for systems with stable Lyapunov exponents.\par 
	
	In this article, we take a closer look at strongly fast invertible systems with bounded coefficients. By linking the dimensions at which a system admits strong fast invertibility to the multiplicities of Lyapunov exponents, we are able to give a full characterization of regular strongly fast invertible systems similar to that of systems with stable Lyapunov exponents. In particular, we show that the stability of Lyapunov exponents implies strong fast invertibility (even in the absence of regularity). Central to our arguments are certain induced systems on spaces of exterior products that represent the evolution of volumes.\par 
	
	Finally, we derive convergence results for the computation of Lyapunov exponents via Benettin's algorithm using perturbation theory. While the stronger assumption of stable Lyapunov exponents clearly leaves more freedom on how to choose stepsizes, we derive conditions for the stepsizes with which convergence can be ensured even if a system is only strongly fast invertible.\par 
	
\end{abstract}
\mbox{}\\





\tableofcontents

\hypersetup{linkcolor=red}

	\section{Introduction}

	\emph{Strong fast invertibility} is a new property for dynamical systems introduced by Anthony Quas, Philippe Thieullen and Mohamed Zarrabi in 2019 \cite{QuasEtAl2019ExplicitBoundsSeparation}. While its original purpose was to ensure uniform invertibility of a cocycle along its fastest growing direction, in our eyes the true benefit lies in its implications for volumes: $L$-dimensional strong fast invertibility implies that the maximal growth of volumes of dimension $L$ from time $s$ to $t$ is (up to a constant) the same as the maximal growth of volumes measured with respect to any initial time $\tau$, i.e., the maximal growth of volumes from time $\tau$ to $t$ divided by the maximal growth of volumes from time $\tau$ to $s$.\par 
	
	In combination with a uniform singular value gap, Quas et al.\ prove the existence of a \emph{dominated equivariant uniform splitting} of the dynamics, i.e., a splitting into equivariant fast and slow subspaces such that the angle between them is bounded from below and solutions corresponding to the fast subspace grow uniformly exponentially faster than solutions corresponding to the slow subspace. Using an equivalent notion, one may also call the splitting \emph{integrally separated} \cite[Def.~2.3]{Sergeev1986ContributionTheoryLyapunov}.\par 
	
	The latter terminology is usually applied in studies of Lyapunov exponents. In particular, the existence of an integrally separated splitting is necessary to ensure the stability of Lyapunov exponents. Essentially, given an integrally separated splitting, the stability of Lyapunov exponents boils down to the stability on the subspaces of the splitting:
	\begin{theorem*}[{\cite[Thm.~5.4.9]{Adrianova1995IntroductionLinearSystems},\cite{BarabanovDenisenko2007NecessarySufficientConditions}}]
		Assume a linear system with bounded coefficients and Lyapunov exponents $\lambda_1>\dots>\lambda_p$ with multiplici{\-}ties $d_1+\dots+d_p=d$. The Lyapunov exponents are stable if and only if there exists a Lyapunov transformation reducing the system to block diagonal form
		\begin{equation*}
			\dot{y}=\textnormal{diag}(B_1(t),\dots,B_p(t))y,
		\end{equation*}
		where $B_i(t)\in\mathbb{R}^{d_i\times d_i}$ is upper triangular, such that the following hold:
		\begin{enumerate}[label=(\roman*)]
			\item all non-trivial solutions of $\dot{y}_i=B_i(t) y_i$ have characteristic exponent $\lambda_i$,
			\item $\lambda_i$ is stable for $\dot{y}_i=B_i(t) y_i$,
			\item there are constants $a,b>0$ such that
			\begin{equation*}
				\|Y_i(t,s)^{-1}\|^{-1}\geq b e^{a(t-s)}\|Y_{i+1}(t,s)\|
			\end{equation*}
			for all $t\geq s$, where $Y_i(t,s)$ denotes the Cauchy matrix of $\dot{y}_i=B_i(t) y_i$.
		\end{enumerate}
	\end{theorem*}
	
	One of our two main goals is to provide a characterization of strong fast invertibility that allows a direct comparison to stability of Lyapunov exponents via the above theorem. To connect strong fast invertibility, a property that concerns the evolution of volumes, to Lyapunov exponents, we work with certain induced systems on spaces of exterior products. These systems allow us to link Lyapunov exponents to volume growth assuming the original system is regular. Our characterization theorem is the following:
	
	\begin{theorem*}
		Assume a linear system with bounded coefficients and Lyapunov exponents $\lambda_1>\dots>\lambda_p$ with multiplicities $d_1+\dots+d_p=d$. If the system is regular and strongly fast invertible at dimensions $d_1+\dots+d_l$ for $l=1,\dots,p$, then there exists a Lyapunov transformation reducing the system to block diagonal form
		\begin{equation*}
			\dot{y}=\textnormal{diag}(B_1(t),\dots,B_p(t))y,
		\end{equation*}
		where $B_i(t)\in\mathbb{R}^{d_i\times d_i}$ is upper triangular, such that the following hold:
		\begin{enumerate}[label=(\roman*)]
			\item all non-trivial solutions of $\dot{y}_i=B_i(t) y_i$ have characteristic exponent $\lambda_i$,
			\item there is a constant $b>0$ such that
			\begin{equation*}
				\|Y_i(t,s)^{-1}\|^{-1}\geq b\|Y_{i+1}(t,s)\|
			\end{equation*}
			for all $t\geq s$, where $Y_i(t,s)$ denotes the Cauchy matrix of $\dot{y}_i=B_i(t) y_i$.
		\end{enumerate}
		Conversely, any block diagonal system $\dot{y}=B(t)y$ satisfying $(i)$ and $(ii)$ is strongly fast invertible at dimensions $d_1+\dots+d_l$ for $l=1,\dots,p$.
	\end{theorem*}  
	
	In particular, our theorem provides an equivalent characterization of strong fast invertibility for systems that are regular and have bounded coefficients. Moreover, it shows that the stability of Lyapunov exponents implies strong fast invertibility at the respective dimensions. More aspects, such as another characterization in the case of simple Lyapunov spectra, can be found in our article. However, we note that there are still interesting aspects to explore that we did not pursue here.\par 
	
	Our second main goal is to derive convergence results for the computation of Lyapunov exponents. More precisely, we focus on \emph{Benettin's algorithm} \cite{BenettinEtAl1980LyapunovCharacteristicExponents, BenettinEtAl1980LyapunovCharacteristicExponentsa} as it is the most fundamental and common algorithm to compute Lyapunov exponents. Its underlying idea is to propagate a set of linear perturbations that are reorthonormalized periodically. The Lyapunov exponents are then computed as averages of volume expansion via the rescaling factors from the orthonormalization procedure.\par 
	
	While it is not difficult to prove convergence of Benettin's algorithm in the absence of numerical errors, integration errors can accumulate and persist. This happens especially if the stepsizes are kept constant. In practice, it is hard to quantify these error, since the exact Lyapunov exponents are usually unknown.\par 
	
	Major efforts have been made by Dieci, Van Vleck and co-authors. They advocate the use of adaptive stepsizes to bound the local integration error and were able to prove error estimates for computed Lyapunov exponents of linear systems that are regular and have stable Lyapunov exponents \cite{DiecivanVleck2005ErrorComputingLyapunov,DiecivanVleck2006PerturbationTheoryApproximation}. While Dieci and Van Vleck proved that the asymptotic limits of the computed exponents can be made arbitrarily close to the true Lyapunov exponents by decreasing the error tolerance or the fixed stepsize, true convergence requires to simultaneously increase the integration time and decrease the stepsizes. This was already conjectured by Mc Donald and Higham in their error analysis for autonomous linear systems in 2001 \cite[Sec.~5]{McdonaldHigham2001ErrorAnalysisQR}.\par 
	
	By tackling both limits simultaneously, we derive new convergence results that differentiate between systems with stable Lyapunov exponents and systems that are only strongly fast invertible. The main difference between the respective convergence results are the requirements for stepsizes. While stepsizes $h_n$ such that
	$$\sum_{n=1}^{\infty}h_n=\infty\quad\text{and}\quad h_n\to 0$$
	are enough to achieve convergence for systems with stable Lyapunov exponents, we need stricter assumptions to compute Lyapunov exponents for strongly fast invertible systems:
	$$\sum_{n=1}^{\infty}h_n=\infty\quad\text{and}\quad \sum_{n=1}^{\infty}h_n^{p+1}<\infty,$$
	where $p>0$ is the order of consistency of the numerical integrator.\par

\bigskip

	\section{Exterior products and powers}

	Exterior products are a handy tool when it comes to studying the evolution of volumes and hence also Lyapunov exponents. In this section we briefly introduce them and some of their properties. Our main reference is \cite[Sec.~3.2.3]{Arnold1998RandomDynamicalSystems}.\par
	
	For $1\leq L\leq d$, the $L$-fold \emph{exterior power} of $\mathbb{R}^d$ is the space $\wedge^L\mathbb{R}^d$ consisting of alternating $L$-linear forms on the dual space $(\mathbb{R}^d)^*\cong \mathbb{R}^d$. A basis can be obtained by taking \emph{exterior products} of basis elements of $\mathbb{R}^d$. For example, the set
	\begin{align*}
		&\{e_I:=e_{i_1}\wedge\dots\wedge e_{i_L}\ |\ I=(i_1,\dots,i_L)\text{ with}\\
		&\hspace{2em}1\leq i_1<\dots<i_L\leq d\},
	\end{align*}
	where $e_i$ is the $i$-th unit vector of $\mathbb{R}^d$, defines a natural basis of $\wedge^L\mathbb{R}^d$. In particular, $\wedge^L\mathbb{R}^d$ has dimension $\binom{d}{L}$.\par 
	
	Not all elements of $\wedge^L\mathbb{R}^d$ are \emph{decomposable}, i.e., of the form $u_1\wedge\dots\wedge u_L$. Some elements are \emph{indecomposable} and can only be expressed as linear combinations of decomposable elements.\par 
	
	Given subspaces $U,V\subset \mathbb{R}^d$, we define induced subspaces
	\begin{align*}
		&(\wedge^kU)\wedge(\wedge^{L-k}V):=\textnormal{span}\{u_1\wedge\dots\wedge u_k \\
		&\hspace{2em}\wedge v_{k+1}\wedge\dots\wedge v_L\ |\ u_i\in U, v_i\in V\}\subset\wedge^L\mathbb{R}^d
	\end{align*}
	as spans of the corresponding induced decomposable elements.\par
	
	By bilinear extension from the set of decomposable elements to $\wedge^L\mathbb{R}^d$, the following defines a scalar product on $\wedge^L\mathbb{R}^d$:
	\begin{equation*}
		\langle u_1\wedge\dots\wedge u_L,v_1\wedge\dots\wedge v_L \rangle:=\det (\langle u_i,v_j\rangle)_{ij}.
	\end{equation*} 
	In particular, the induced norm of a decomposable element $u_1\wedge\dots\wedge u_L$ is the $L$-volume of the parallelepiped spanned by $u_1,\dots,u_L$:
	\begin{equation*}
		\|u_1\wedge\dots\wedge u_L\|=\sqrt{\det (\langle u_i,u_j\rangle)_{ij}}.
	\end{equation*}
	
	Throughout this article norms without annotation always denote euclidean norms or associated matrix norms.\par 
	
	\newpage
	
	\begin{lemma}\label[lemma]{lemma:ExteriorPowerVectorProperties}	
		We have 
		\begin{enumerate}
			\item[(i)] $\|u_1\wedge\dots\wedge u_L\|\leq \|u_1\wedge\dots\wedge u_k\|$\\
			$\cdot\|u_{k+1}\wedge\dots\wedge u_L\|$
		\end{enumerate}
		for $u_1,\dots,u_L\in\mathbb{R}^d$ and
		\begin{enumerate}
			\item[(ii)] $\langle \hat{u}\wedge\hat{v},\hat{u}'\wedge\hat{v}'\rangle=\langle\hat{u},\hat{u}'\rangle\langle\hat{v},\hat{v}'\rangle$,
			\item[(iii)] $\|\hat{u}\wedge\hat{v}\|= \|\hat{u}\|\,\|\hat{v}\|$
		\end{enumerate}
		for $\hat{u},\hat{u}'\in \wedge^k U$ and $\hat{v},\hat{v}'\in\wedge^{L-k}V$ with $U,V\subset\mathbb{R}^d$ orthogonal.
	\end{lemma}
	
	\begin{proof}
		$(i)$ can be found in \cite[Subsec.~3.2.3]{Arnold1998RandomDynamicalSystems}.\par 
		
		Let $U,V\subset\mathbb{R}^d$ with $U\perp V$. Choose bases $(u_i)_i$ of $U$ and $(v_j)_j$ of $V$. We denote the elements of the induced bases of $\wedge^kU$ and $\wedge^{L-k}V$ by $u_I$ and $v_J$ respectively. Since
		\begin{equation*}
			\langle u_I\wedge v_J, u_{I'}\wedge v_{J'}\rangle
		\end{equation*}
		is the determinant of a block diagonal matrix, it is equal to the product of the determinants of both blocks, which is
		\begin{equation*}
			\langle u_I,u_{I'}\rangle\langle v_J,v_{J'}\rangle.
		\end{equation*}
		Thus, $(ii)$ holds on basis elements and by bilinearity of the inner product on arbitrary elements.\par 
		
		Assertion $(iii)$ follows immediately from $(ii)$ by setting $\hat{u}=\hat{u}'$ and $\hat{v}=\hat{v}'$. 
	\end{proof}
	
	Next, we discuss several constructions for linear maps. The $L$-fold exterior power of $A\in\mathbb{R}^{d\times d}$ is defined via 
	\begin{equation*}
		(\wedge^LA)(u_1\wedge\dots\wedge u_L):=Au_1\wedge\dots \wedge Au_L.
	\end{equation*}
	Similarly, one may define
	\begin{align*}
		&((\wedge^kA)\wedge(\wedge^{L-k}B))(u_1\wedge\dots\wedge u_L)\\
		&\hspace{1em}:=Au_1\wedge\dots\wedge Au_k\wedge Bu_{k+1}\wedge\dots \wedge Bu_L.
	\end{align*}
	for $A,B\in\mathbb{R}^{d\times d}$. Another helpful construction is
	\begin{align*}
		&\hat{A}^L(u_1\wedge\dots\wedge u_L):=\\
		&\hspace{1em}\sum_{k=1}^{L}u_1\wedge\dots\wedge u_{k-1}\wedge Au_k\wedge u_{k+1}\wedge\dots\wedge u_L.
	\end{align*}
	
	Since we will make extensive use of these induced maps, we state and derive a list of basic properties for them.\par
	
	\begin{lemma}\label[lemma]{lemma:ExteriorPowerMatrixProperties}	
		The following are true for $A,B\in\mathbb{R}^{d\times d}$:
		\begin{enumerate}[label=(\roman*)]
			\item $\wedge^L I_{\mathbb{R}^d} = I_{\wedge^L\mathbb{R}^d}$,
			\item $\wedge^L(AB)=(\wedge^LA)(\wedge^LB)$,
			\item $(\wedge^LA)^{-1}=\wedge^LA^{-1}$ if $A\in GL(d,\mathbb{R})$,
			\item $\reallywidehat{\alpha A + \beta B}^L=\alpha\hat{A}^L+\beta\hat{B}^L$,
			\item $\|\wedge^LA\|=\sigma_1(A)\dots\sigma_L(A)$, where \\$\sigma_1\geq\dots\geq\sigma_d$ denote the singular values,
			\item $\|\wedge^LA\|\leq \|A\|^L$,
			\item $\|\hat{A}^L\|\leq L\|A\|$,
			
			\item \small $\|(\wedge^kA)\wedge (\wedge^{L-k}B)\|\leq \binom{d}{L}^{\frac{1}{2}}\|\wedge^kA\|\,\|\wedge^{L-k}B\|$,\normalsize
			\item if $A[u_1,\dots,u_L]=QR$, then\\
			$\|(\wedge^LA)(u_1\wedge\dots\wedge u_L)\|=r_{11}\dots r_{LL}$, where $r_{ii}\geq 0$ denote the diagonal elements of $R$ ordered in decreasing size,
			\item $\det(\wedge^LA)= \det(A)^{\binom{d-1}{L-1}}$.
		\end{enumerate}		
	\end{lemma}
	
	\begin{proof}
		$(i)-(vii)$ and $(x)$ can be found in \cite[Subsec.~3.2.3]{Arnold1998RandomDynamicalSystems}.\par 
		
		To show $(viii)$, write
		\begin{equation*}
			I_1 = (i_1,\dots,i_k)\text{ and }I_2=(i_{k+1},\dots,i_L)
		\end{equation*}
		for a given tuple $I=(i_1,\dots,i_L)$ and estimate
		\begin{align*}
			&\|(\wedge^kA)\wedge (\wedge^{L-k}B)\|\\
			&\hspace{1em}=\sup_{\sum \alpha_I^2=1}\left\|((\wedge^kA)\wedge (\wedge^{L-k}B))\left(\sum\alpha_Ie_I\right)\right\|\\
			&\hspace{1em}\leq \sup_{\sum \alpha_I^2=1}\sum|\alpha_I|\,\|((\wedge^kA)e_{I_1})\wedge((\wedge^{L-k}B)e_{I_2})\|\\
			&\hspace{1em}\leq \sup_{\sum \alpha_I^2=1}\sum|\alpha_I|\,\|(\wedge^kA)e_{I_1}\|\,\|(\wedge^{L-k}B)e_{I_2}\|\\
			&\hspace{1em}\leq \binom{d}{L}^{\frac{1}{2}}\|\wedge^kA\|\,\|\wedge^{L-k}B\|.
		\end{align*}
		Finally, we prove $(ix)$:
		\begin{align*}
			&\|(\wedge^LA)(u_1\wedge\dots\wedge u_L)\|^2=\det(\langle Au_i,Au_j\rangle)_{ij}\\
			&\hspace{1em}=\det((QR)^T(QR))=\det(R^TR)=r_{11}^2\dots r_{LL}^2.
		\end{align*}	
	\end{proof}
	
	\begin{lemma}\label[lemma]{lemma:ExtPowerPropertiesDiag}
		Let $A=\textnormal{diag}(A_1,A_2)\in\mathbb{R}^{d\times d}$ and set $B:=\textnormal{diag}(A_1,0)$ and $C:=\textnormal{diag}(0,A_2)$. The following are true:
		\begin{enumerate}[label=(\roman*)]
			\item $\wedge^LA=\sum_{k=0}^{L}(\wedge^kB)\wedge(\wedge^{L-k}C)$,
			\item $\|\wedge^LA\|=\max_k\|(\wedge^kB)\wedge(\wedge^{L-k}C)\|$,
			\item $\|(\wedge^kB)\wedge(\wedge^{L-k}C)\|\geq \|\wedge^kB\|\,\|\wedge^{L-k}C\|$.
		\end{enumerate}
	\end{lemma}
	
	\begin{proof}
		
		Let $A_1\in\mathbb{R}^{d_1\times d_1}$. Given a basis element
		\begin{equation*}
			e_I=e_{i_1}\wedge\dots\wedge e_{i_L},
		\end{equation*}
		choose $j$ such that $i_j\leq d_1<i_{j+1}$. It holds
		\begin{align*}
			(\wedge^LA)e_I=((\wedge^jB)\wedge(\wedge^{L-j}C))e_I.
		\end{align*}
		Moreover, all summands in 
		\begin{equation*}
			\sum_{k=0}^{L}((\wedge^kB)\wedge(\wedge^{L-k}C))e_I
		\end{equation*}
		vanish except for $k=j$. Hence, the two maps in $(i)$ coincide on basis elements $e_I$.\par 
		
		The subspaces
		\begin{equation*}
			W_k:=\left(\wedge^k\left(\mathbb{R}^{d_1}\times\{0\}\right)\right)\wedge\left(\wedge^{L-k}\left(\{0\}\times\mathbb{R}^{d-d_1}\right)\right)
		\end{equation*}
		for $k=0,\dots,L$ form an orthogonal decomposition of $\wedge^L\mathbb{R}^d$. Since
		$$\textnormal{im}\left((\wedge^kB)\wedge(\wedge^{L-k}C)\right)\subset W_k$$
		and 
		$$W_i\subset \ker\left((\wedge^kB)\wedge(\wedge^{L-k}C)\right)$$
		for $i\neq k$, $(ii)$ easily follows by means of this decomposition.\par 
		
		To prove $(iii)$, we use $(iii)$ of \cref{lemma:ExteriorPowerVectorProperties}. It follows that
		\begin{align*}
			&\|(\wedge^kA)\wedge (\wedge^{L-k}B)\|\\
			&\hspace{1em}\geq\max_{\|\hat{u}\wedge \hat{v}\|=1}\|(\wedge^kA)\hat{u}\wedge(\wedge^{L-k}B)\hat{v}\|\\
			&\hspace{1em}=\max_{\|\hat{u}\|=1}\|(\wedge^kA)\hat{u}\|\max_{\|\hat{v}\|=1}\|(\wedge^{L-k}B)\hat{v}\|\\
			&\hspace{1em}=\|\wedge^kA\|\,\|\wedge^{L-k}B\|,
		\end{align*}
		where maxima are with respect to $\hat{u}\wedge\hat{v}\in W_k$.
	\end{proof}
	
	Next, we relate the principle angles between two complementary subspaces of $\mathbb{R}^d$ to a principle angle on $\wedge^L\mathbb{R}^d$.
	
	\begin{proposition}\label[proposition]{prop:PrincipleAngles}
		Let $\mathbb{R}^d=U\oplus V$ with $\dim U=L$. Set $\hat{U}:=\wedge^LU$ and $\hat{V}:=\left(\wedge^LV^{\perp}\right)^{\perp}$. Then
		\begin{equation}\label{eqn:PrincipleAngles}
			\prod_{i=1}^{\min(L,d-L)}\sin\alpha_i(U,V)=\sin \alpha_1(\hat{U},\hat{V}),
		\end{equation}
		where $0<\alpha_1\leq\alpha_2\leq \dots\leq  \frac{\pi}{2}$ are the principle angles between the respective subspaces.
	\end{proposition}
	
	\begin{proof}
		First we pass to the orthogonal complement of $V$ in order to work with two subspaces of the same dimension. According to \cite{ZhuKnyazev2013AnglesSubspacesTheir} it holds
		\begin{align*}
			\prod_{i=1}^{\min(L,d-L)}\sin \alpha_i(U,V)&=\prod_{i=1}^{L}\sin\left(\frac{\pi}{2}-\alpha_i(U,V^{\perp})\right)\\
			&=\prod_{i=1}^{L}\cos\alpha_i(U,V^{\perp}).
		\end{align*}
		The concept of ``higher dimensional angle'' coined in \cite{Jiang1996AnglesEuclideanSubspaces} helps us to transition to $\wedge^L\mathbb{R}^d$. In fact, the higher dimensional angle $\theta(U,V^{\perp})$ between two subspaces $U$ and $V^{\perp}$ of the same dimension satisfies
		\begin{equation*}
			\cos\theta(U,V^{\perp})=\prod_{i=1}^{L}\cos\alpha_i(U,V^{\perp})
		\end{equation*}
		and is defined via 
		\begin{equation*}
			\cos\theta(U,V^{\perp}):=\frac{\langle \alpha,\beta\rangle}{\|\alpha\|\,\|\beta\|},
		\end{equation*}
		where $\alpha,\beta\in\wedge^L\mathbb{R}^d$ are decomposable elements representing $U$ and $V^{\perp}$, i.e., $\alpha=u_1\wedge\dots\wedge u_L$ with $\textnormal{span}(u_1,\dots,u_L)=U$ and $\beta=v_1'\wedge\dots\wedge v_L'$ with $\textnormal{span}(v_1',\dots,v_L')=V^{\perp}$. The latter is nothing else than the angle between $\alpha$ and $\beta$ or the principle angle between their corresponding $1$-dimensional subspaces:
		\begin{align*}
			\cos\theta(U,V^{\perp})&=\cos\alpha_1(\wedge^LU,\wedge^LV^{\perp})\\
			&=\cos\alpha_1(\hat{U},\hat{V}^{\perp})\\
			&=\sin\alpha_1(\hat{U},\hat{V}).
		\end{align*}	
	\end{proof}
	
	\begin{corollary}\label[corollary]{cor:PrincipleAngles}
		In \cref{prop:PrincipleAngles} it holds
		\begin{equation*}
			\left(\wedge^LV^{\perp}\right)^{\perp}=\oplus_{k=0}^{L-1}\left((\wedge^{k}U)\wedge(\wedge^{L-k}V)\right).
		\end{equation*}
		In particular, \cref{eqn:PrincipleAngles} is true for
		\begin{align*}
			U&:=\textnormal{span}\{Ae_1,\dots,Ae_L\},\\
			V&:=\textnormal{span}\{Ae_{L+1},\dots,Ae_d\},\\
			\hat{U}&=\textnormal{span}\{(\wedge^LA)(e_1\wedge\dots\wedge e_L)\},\\
			\hat{V}&=\textnormal{span}\{(\wedge^LA)(e_{i_1}\wedge\dots\wedge e_{i_L})\ |\\
			&\hspace{2em}\{i_1,\dots,i_L\}\neq \{1,\dots,L\}\}
		\end{align*}
		for any given $A\in GL(d,\mathbb{R})$.
	\end{corollary}
	
	\begin{proof}
		The identity can be checked by using that both spaces have codimension $1$ and by showing that the inner product between elements of the form 
		$$v_1'\wedge\dots\wedge v_L'\in\wedge^LV^{\perp}$$
		and
		$$u_1\wedge\dots\wedge u_k\wedge v_{k+1}\wedge\dots\wedge v_L$$
		with $u_i\in U$ and $v_i\in V$ vanishes for $k<L$.\par 
	\end{proof}

\bigskip

	\section{Asymptotic properties of li{\-}near systems}

	Consider the linear differential equation 
	\begin{equation}\label{eqn:system}
		\dot{x}=A(t)x, \quad A\in C(\mathbb{R}_{\geq 0},\mathbb{R}^{d\times d})
	\end{equation}
	with bounded coefficients, i.e., 
	$$M:=\sup_t \|A(t)\|<\infty.$$
	
	We write $X(t)$ for a fundamental matrix of \cref{eqn:system} and $X(t,s):=X(t)X(s)^{-1}$, $t\geq s$, for the Cauchy matrix. By requiring that $A$ is bounded, solutions grow or decay at most exponentially fast:
	\begin{equation*}
		e^{-M(t-s)}\leq\|X(t,s)^{\pm 1}\|\leq e^{M(t-s)}
	\end{equation*}
	for all $t\geq s$.\par 
	
	In our article, \cref{eqn:system} will be changed by so-called \emph{Lyapunov transformations}. These transformations leave asymptotic properties\footnote{We will define all mentioned properties throughout our article.} like boundedness of the system matrix, the Lyapunov spectrum, regularity, the stability of Lyapunov exponents, or strong fast invertibility invariant.\par 
	
	\begin{definition}[{\cite[Def.~3.1.1]{Adrianova1995IntroductionLinearSystems}}]
		$x=L(t)y$ is called \emph{Lyapunov transformation} if $L$ is continuously differentiable, everywhere invertible and $L,L^{-1},\dot{L}$ are bounded. The transformed system $\dot{y}=B(t)y$ is given by
		\begin{equation*}
			B(t)=L(t)^{-1}A(t)L(t)-L(t)^{-1}\dot{L}(t).
		\end{equation*}
	\end{definition}
	
	It is easy to see that Lyapunov transformations form a group with respect to multiplication. A particular Lyapunov transformation is obtained via the Gram-Schmidt procedure.\par
	
	\begin{proposition}[{\small\cite[Lemma~3.3.1 and Thm.~3.3.1]{Adrianova1995IntroductionLinearSystems}}]\label[proposition]{proposition:GramSchmidt}
		There is an orthogonal Lyapunov transformation such that the fundamental matrix $Y(t)$ of the transformed system is upper triangular.
	\end{proposition}

\subsection{Characteristic exponents}
	
	In this subsection, we define characteristic exponents and Lyapunov exponents.\par 
	
	\begin{definition}[{\cite[Def.~2.1.1]{Adrianova1995IntroductionLinearSystems}}]
		The \emph{characteristic exponent} of $f:[0,\infty)\to \mathbb{R}$ is given by
		\begin{equation*}
			\chi[f]:=\limsup_{t\to\infty}\frac{1}{t}\log |f(t)|\in\mathbb{R}\cup\{\pm \infty\}.
		\end{equation*}		
	\end{definition}
	
	A few handy properties follow easily from the definition (see \cite[Sec.~2.1]{Adrianova1995IntroductionLinearSystems} and \cite[Lemma~3.2.1]{Arnold1998RandomDynamicalSystems}):
	\begin{enumerate}[label=(\roman*)]
		\item $\chi[c]=0$ for $c\neq 0$
		\item $\chi[cf]=\chi[f]$ for $c\neq 0$,
		\item $\chi[|f|^c]=c\chi[f]$ for $c\in\mathbb{R}$ (set $0(\pm\infty)=0$),
		\item $\chi[f]\leq\chi[g]$ if $|f|\leq |g|$,
		\item $\chi[f+g]\leq \max(\chi[f],\chi[g])$ with equality if $\chi[f]\neq\chi[g]$,
		\item $\chi[fg]\leq \chi[f]+\chi[g]$ (if the right-hand side makes sense).
	\end{enumerate}\par 
	
	In the context of dynamical systems, we apply characteristic exponents to measure the (upper) exponential growth of solutions $x(t)$:
	\begin{equation*}
		\chi[x]:=\chi[\|x\|]=\limsup_{t\to\infty}\frac{1}{t}\log \|x(t)\|.
	\end{equation*}
	Note that $\chi[x]$ is independent of the chosen norm because all norms on $\mathbb{R}^d$ are equivalent.\footnote{In infinite dimensions the characteristic exponents generally depend on the chosen norm, although independence can be achieved on certain scales of Banach spaces \cite{BlumenthalPunshon-Smith2023NormEquivalenceLyapunov}.}\par 
	
	Since solutions with distinct characteristic exponents are linearly independent, the characteristic exponents can take at most $d$ values
	\begin{equation*}
		\infty>M\geq \lambda_1>\dots>\lambda_p\geq -M>-\infty.
	\end{equation*}
	In particular, there is a filtration of subspaces of the space of solutions called \emph{Lyapunov filtration}
	\begin{equation*}
		\mathbb{R}^d\cong V_1\supset \dots\supset V_p\supset  V_{p+1}:=\{0\} 
	\end{equation*}
	with $V_i:=\{x\, |\, \chi[x]\leq\lambda_i\}$ satisfying
	\begin{equation*}
		\chi[x]=\lambda_i\quad\iff\quad x\in V_i\setminus V_{i+1}.
	\end{equation*}
	We set $d_i:=\dim V_i - \dim V_{i+1}$.\par 
	
	\begin{remark}
		Sometimes the Lyapunov filtration is defined on the space of initial vectors. In fact, both versions of filtration spaces are naturally isomorphic:
		$$V_{i}':=\{v\ |\ \chi[X(t,0)v]\leq \lambda_i\}=\{x(0)\ |\ x\in V_i\}.$$
		One may also define the filtration spaces using a given fundamental matrix $X(t)$:
		$$V_{i,X}:=\{v\ |\ \chi[X(t)v]\leq \lambda_i\}.$$
	\end{remark}
	
	\begin{definition}
		The \emph{(forward) Lyapunov spectrum} $(\lambda_i,d_i)_{i=1,\dots,p}$ of a system consists of its \emph{Lyapunov exponents} $\lambda_i$ together with their \emph{multiplicities} $d_i$. Since we sometimes count Lyapunov exponents according to their multiplicities, we define $\Lambda_1,\dots,\Lambda_d$ via
		\begin{equation*}
			\Lambda_{d_1+\dots+d_{i-1}+j}:=\lambda_i\text{ for }j=1,\dots,d_i.
		\end{equation*}
		
		A Lyapunov exponent $\lambda_i$ is called \emph{simple} or \emph{nondegenerate} if $d_i=1$ and otherwise \emph{degenerate}. If all Lyapunov exponents are simple, we call the Lyapunov spectrum \emph{simple}.
	\end{definition}
	
	We may choose\footnote{Construct a basis by iteratively choosing $d_i$ solutions in $V_i\setminus V_{i+1}$ for $i=p,\dots,1$.} a basis of solutions according to the Lyapunov filtration such that it realizes the whole Lyapunov spectrum. In that case
	\begin{equation}\label{eqn:normal}
		\sum_{i=1}^{d}\chi[x_i]=\sum_{i=1}^{p}d_i\lambda_i.
	\end{equation}\par 
	
	\begin{definition}[{\cite[Def.~2.4.2]{Adrianova1995IntroductionLinearSystems}}]
		A basis $x_1,\dots,x_d$ of solutions is called \emph{normal} if $\sum_i\chi[x_i]$ is minimal, i.e., if \cref{eqn:normal} holds. Moreover, we call a fundamental matrix \emph{normal} if its columns form a normal basis.
	\end{definition}
	
	If we order a normal basis such that $\chi[x_i]$ decreases with $i$, then $\chi[x_i]=\Lambda_i$.\par
	
	\begin{proposition}[{\cite[Thm.~2.5.1, Cor.~2.5.1 and Remark~2.5.2]{Adrianova1995IntroductionLinearSystems}}]
		Any fundamental matrix $X(t)=[x_1,\dots,x_d]$ satisfies
		\begin{equation*}
			\sum_{i=1}^d\chi[x_i]\geq \sum_{i=1}^{p}d_i\lambda_i\geq \chi[\det X(t)].
		\end{equation*}
		If equality holds, then $X(t)$ is normal. However, $X(t)$ being normal does not imply equality in general.
	\end{proposition}
	
	The largest Lyapunov exponent can be expressed as the characteristic exponent of $\|X(t)\|$.
	
	\begin{proposition}\label[proposition]{prop:CharExpMatrixNorm}
	The largest Lyapunov exponent satisfies
	\begin{equation*}
		\lambda_1=\chi[\|X(t)\|].
	\end{equation*}
	\end{proposition}
	
	\begin{proof}
		One inequality is trivial and the other follows from
		\begin{equation*}
			\|X(t)\|\leq\sqrt{d}\max_i\|X(t)e_i\|_1,
		\end{equation*}
		where $e_i$ is the $i$-th unit vector of $\mathbb{R}^d$.
	\end{proof}

\subsection{Regularity}
	
	We now introduce the notion of regularity, which implies that characteristic exponents of solutions can be obtained as limits instead of limes superiors.\par 
	
	\begin{definition}[{\cite[p.~115]{Arnold1998RandomDynamicalSystems}}]
		We call a system \emph{(forward) regular}\footnote{The definition from Arnold's book is equivalent to \cite[Def.~3.5.1]{Adrianova1995IntroductionLinearSystems}. This can be checked via \cite[Lemma~3.5.1]{Adrianova1995IntroductionLinearSystems} and the Liouville-Ostrogradski formula $$\det X(t) = \det X(t_0)e^{\int_{t_0}^{t}\textnormal{tr}(A(\tau))\,d\tau}.$$} if 
		\begin{equation*}
			\sum_{i=1}^{p}d_i\lambda_i=\liminf_{t\to\infty}\frac{1}{t}\log|\det X(t)|.
		\end{equation*}
	\end{definition}
	
	In case of triangular systems regularity can be checked via the diagonal elements of the system matrix \cite[Cor.~3.8.1]{Adrianova1995IntroductionLinearSystems}. Moreover, one may check regularity of general systems through Perron's regularity test, which compares the Lyapunov spectra of the system and its adjoint system $\dot{y}=-A(t)^Ty$ \cite[Thm.~3.6.1]{Adrianova1995IntroductionLinearSystems}.\par
	
	\begin{proposition}[{\cite[Thm.~3.9.1]{Adrianova1995IntroductionLinearSystems}}]
		Regular systems have sharp characteristic exponents, i.e., it holds
		\begin{equation*}
			\chi[x]=\lim_{t\to\infty}\frac{1}{t}\log\|x(t)\|.
		\end{equation*}
	\end{proposition}
	
	All Lyapunov exponents of a regular system can be obtained via characteristic exponents of singular values of $X(t)$.\par 
	
	\begin{proposition}\label[proposition]{prop:SingVals}
		If \cref{eqn:system} is regular, then
		\begin{equation*}
			\Lambda_i=\lim_{t\to\infty}\frac{1}{t}\log\sigma_i(X(t))
		\end{equation*}
		for all $i=1,\dots,d$, where $\sigma_1\geq\dots\geq\sigma_d$ denote the singular values. 
	\end{proposition}
	
	We prove \cref{prop:SingVals} in the next subsection using induced systems on the space of exterior products. Next, let us state a result on singular vectors of $X(t)$.\par 
	
	\begin{proposition}\label[proposition]{prop:MET}
		If \cref{eqn:system} is regular and $u_1(t),\dots,u_d(t)$ denote the right singular vectors of $X(t)$, then $$U_{i,X}(t):=\textnormal{span}(u_{d_1+\dots+d_{i-1}+1}(t),\dots,u_{d_{1}+\dots+d_i}(t))$$
		converges to $U_{i,X}:=V_{i+1,X}^{\perp}\cap V_{i,X}$ exponentially fast. More precisely, it holds
		\begin{equation*}
			\limsup_{t\to\infty}\frac{1}{t}\log\max_{\substack{u\in U_{i,X}(t),\\ u'\in U_{j,X},\\\|u\|=\|u'\|=1}}|\langle u,u'\rangle|\leq -|\lambda_i-\lambda_j|
		\end{equation*} 
		for $i\neq j$.
	\end{proposition}
	
	\begin{proof}
		Since \cref{eqn:system} has bounded coefficients and its induced systems are regular (see \cref{lemma:InducedNormalBasis}), the deterministic version of the multiplicative ergodic theorem \cite[Prop.~3.4.2]{Arnold1998RandomDynamicalSystems} applies. In its proof, Arnold shows that the filtration $F_X(t)$ given by the spaces
		$$V_{i,X}(t):=U_{p,X}(t)\oplus\dots\oplus U_{i,X}(t)$$
		converges exponentially fast to the Lyapunov filtration $F_X$ given by the spaces
		$$V_{i,X}=U_{p,X}\oplus\dots\oplus U_{i,X}$$
		using the metric $\delta$ (see \cite[Eq.~(3.4.10)]{Arnold1998RandomDynamicalSystems}) on the manifold of filtrations. In particular, he shows that
		$$\limsup_{n\to\infty}\frac{1}{n}\log \delta(F_X(n),F_X)\leq -h,$$
		where $h>0$ is a parameter also appearing in the definition of $\delta$. Disentangling the metric yields the desired convergence statement for discrete time.\par 
		The version for continuous time follows as described in the proof of \cite[Thm.~3.4.1]{Arnold1998RandomDynamicalSystems} since 
		$$\limsup_{n\to\infty}\frac{1}{n}\sup_{0\leq t\leq 1}\log \|X(n+t,n)^{\pm 1}\|\leq 0.$$
	\end{proof}
	
	If at least two singular values coincide, the singular value decomposition is not unique. However, according to \cref{prop:SingVals} the singular values corresponding to different Lyapunov exponents are distinct for large $t$. In particular, the spaces $U_{i,X}(t)$ are uniquely defined for large $t$.\par 
	
	\begin{corollary}\label[corollary]{cor:StrongestGrowth}
		If \cref{eqn:system} is regular and $\lambda_1$ is simple, then
		\begin{equation*}
			\lim_{t\to\infty} \frac{\|X(t)v\|}{\|X(t)\|}=\|v-P_{V_{2,X}}v\|
		\end{equation*}
		for $v\in\mathbb{R}^d$, where $P_{V_{2,X}}$ denotes the orthogonal projection onto $V_{2,X}$. 
	\end{corollary}
	
	\begin{proof}
		To prove the corollary, we use the two orthogonal decompositions $\mathbb{R}^d=U_{1,X}\oplus V_{2,X}$ and $\mathbb{R}^d=U_{1,X}(t)\oplus V_{2,X}(t)$ from \cref{prop:MET} and its proof. Since $d_1=1$, there are unit vectors $u_1$ and $u_1(t)$ spanning $U_{1,X}$ and $U_{1,X}(t)$. Moreover, \cref{prop:MET} implies
		$$\lim_{t\to\infty}\|P_{V_{2,X}(t)}u_1\|=0,$$
		where $P_{V_{2,X}(t)}$ denotes the orthogonal projection onto $V_{2,X}(t)$. We estimate $\|X(t)u_1\|\leq \|X(t)\|$ and
		\begin{align*}
			&\|X(t)u_1\|\\
			&\hspace{1em}=\|X(t)(\langle u_1,u_1(t)\rangle u_1(t) + P_{V_{2,X}(t)}u_1)\|\\
			&\hspace{1em}\geq |\langle u_1,u_1(t)\rangle|\,\|X(t)\| - \|X(t)\|\,\|P_{V_{2,X}(t)}u_1\|\\
			&\hspace{1em}=\left(\sqrt{1-\|P_{V_{2,X}(t)}u_1\|^2}-\|P_{V_{2,X}(t)}u_1\|\right)\|X(t)\|.
		\end{align*}
		Thus, it holds
		$$\lim_{t\to\infty}\frac{\|X(t)u_1\|}{\|X(t)\|}=1.$$
		Now, we decompose $v\in\mathbb{R}^d$ into $v=\langle v,u_1\rangle u_1 + v_2$ according to $\mathbb{R}^d=U_{1,X}\oplus V_{2,X}$. Since the system is regular and the characteristic exponents satisfy $\chi[X(t)v_2]\leq \lambda_2$ and $\chi[\|X(t)\|]=\lambda_1$, we have
		$$\lim_{t\to\infty}\frac{\|X(t)v_2\|}{\|X(t)\|}=0.$$
		The claim of the corollary follows from		
		\begin{align*}
			\frac{\|X(t)v\|}{\|X(t)\|}\leq |\langle v,u_1\rangle|\frac{\|X(t)u_1\|}{\|X(t)\|} + \frac{\|X(t)v_2\|}{\|X(t)\|}
		\end{align*}
		and 
		$$\frac{\|X(t)v\|}{\|X(t)\|}\geq |\langle v,u_1\rangle|\frac{\|X(t)u_1\|}{\|X(t)\|} - \frac{\|X(t)v_2\|}{\|X(t)\|}.$$
	\end{proof}

\subsection{Induced systems}
	
	\cref{eqn:system} induces a system on $\wedge^L\mathbb{R}^d$ via
	\begin{equation}\label{eqn:InducedSystem}
		\dot{\hat{x}}=\widehat{A(t)}^L\hat{x}.
	\end{equation}
	In the following, we refer to \cref{eqn:InducedSystem} using the term \emph{induced system}.\par 
	
	A fundamental matrix of the induced system can be obtained by taking the exterior power of a fundamental matrix of the original system:\par 
	
	\begin{proposition}
		It holds
		\begin{equation*}
			\frac{d}{dt}(\wedge^LX)=\hat{A}^L(\wedge^LX).
		\end{equation*}
	\end{proposition}
	
	\begin{proof}
		Since the determinant is multilinear, we have
		\begin{align*}
			&\left\langle\frac{d}{dt}(\wedge^LX)e_I,e_J\right\rangle\\
			&\hspace{1em}=\frac{d}{dt}\left\langle Xe_{i_1}\wedge\dots\wedge Xe_{i_L},e_J\right\rangle\\
			&\hspace{1em}=\sum_{k=1}^{L}\left\langle Xe_{i_1}\wedge\dots\wedge\frac{d}{dt}(Xe_{i_k})\wedge\dots\wedge Xe_{i_L},e_J\right\rangle\\
			&\hspace{1em}=\left\langle\sum_{k=1}^{L} Xe_{i_1}\wedge\dots\wedge AXe_{i_k}\wedge\dots\wedge Xe_{i_L},e_J\right\rangle\\
			&\hspace{1em}=\left\langle\hat{A}^L(\wedge^LX)e_I,e_J\right\rangle.
		\end{align*}
	\end{proof}
	
	Several properties carry over from the original system to the induced system. For instance, the Lyapunov spectrum $(\lambda_{i,L},d_{i,L})_{i=1,\dots,p_L}$ of the induced system can be related to the spectrum of the original system if the system is regular (see also \cite[Thm.~5.3.1]{Arnold1998RandomDynamicalSystems}).\par 
	
	\begin{theorem}\label{thm:METExtPowers}
		Let \cref{eqn:system} be regular. The Lyapunov exponents $\lambda_{i,L}$ of the induced system are the different values given by 
		\begin{equation*}
			\Lambda_{i_1}+\dots+\Lambda_{i_L}
		\end{equation*}
		for indices $i_1<\dots<i_L$. The corresponding multiplicity $d_{i,L}$ is the number of combinations of indices for which $\lambda_{i,L}$ can be achieved. In particular,
		\begin{equation*}
			\lambda_{1,L} = \Lambda_1 + \dots + \Lambda_L.
		\end{equation*}
		If $L=d_1+\dots+d_l$, then $\lambda_{1,L}$ is simple and the second space of the Lyapunov filtration $V_{2,L}$ is given by 
		\small
		\begin{align*}
			\textnormal{span}\{x_1\wedge\dots\wedge x_L\ |\ \textnormal{span}(x_1,\dots,x_L)\cap V_{l+1}\neq \{0\}\}.
		\end{align*}
		\normalsize
	\end{theorem}
	
	\begin{remark}
		If \cref{eqn:system} is not regular, the spectrum of the induced system can differ from what is described in \cref{thm:METExtPowers}.
	\end{remark} 
	
	\begin{proof}
		In the absence of regularity, the fastest growing solutions do not necessarily span the fastest growing subspace. This can be the case if the characteristic exponents of the fastest growing solutions are only obtainable along distinct subsequences.\par 
		
		Indeed, one readily checks that
		\begin{align*}
			X(t):=\begin{pmatrix}
				e^{t\sin(\log(1+t))} & 0 & 0\\
				0 & e^{t\cos(\log(1+t))} & 0\\
				0 & 0 & e^{\frac{1}{2}t}
			\end{pmatrix}
		\end{align*}
		is a normal fundamental matrix corresponding to a bounded, continuous system with Lyapunov exponents $\Lambda_1=\Lambda_2=1$ and $\Lambda_3=1/2$. Since the characteristic exponents of the first two columns of $X=[x_1,x_2,x_3]$ are not sharp, the system is not regular. Moreover, one computes
		\begin{align*}
			\chi[x_1\wedge x_2]&=\chi[e^{t(\sin(\log(1+t))+\cos(\log(1+t)))}]=\sqrt{2}\\
			\chi[x_1\wedge x_3]&=\chi[e^{t(\sin(\log(1+t))+\frac{1}{2})}]=\frac{3}{2}\\
			\chi[x_2\wedge x_3]&=\chi[e^{t(\cos(\log(1+t))+\frac{1}{2})}]=\frac{3}{2}.			
		\end{align*}
		Since $x_1(t)\wedge x_3(t)$ and $x_2(t)\wedge x_3(t)$ are orthogonal for each $t$, any nontrivial linear combination of them has characteristic exponent $3/2$. In particular, $\Lambda_{1,2}=\Lambda_{2,2}=3/2$ and $\Lambda_{3,2}=\sqrt{2}$.	Thus, even though $x_1$ and $x_2$ have the highest characteristic exponents, their associated volume element has the lowest characteristic exponent.
	\end{proof}
	
	We prove \cref{thm:METExtPowers} using the following lemma:\par 
	
	\begin{lemma}\label[lemma]{lemma:InducedNormalBasis}
		If \cref{eqn:system} is regular, then the induced basis of a normal basis is normal and the induced system is regular.
	\end{lemma}
	
	\begin{proof}
		Let $x_1,\dots,x_d$ be a normal basis such that $\chi[x_i]=\Lambda_i$. Since
		\begin{align*}
			\chi[x_{i_1}\wedge\dots\wedge x_{i_L}]&\leq \chi[x_{i_1}]+\dots+\chi[x_{i_L}]\\
			&=\Lambda_{i_1}+\dots+\Lambda_{i_L},
		\end{align*} 
		we have
		\begin{align*}
			&\sum_{i_1<\dots<i_L}\Lambda_{i_1}+\dots+\Lambda_{i_L}\\
			&\hspace{1em}\geq\sum_{i_1<\dots<i_L}\chi[x_{i_1}\wedge\dots\wedge x_{i_L}]\\
			&\hspace{1em}\geq \chi[\det(\wedge^LX(t))]\\
			&\hspace{1em}\geq \liminf_{t\to\infty}\frac{1}{t}\log|\det(\wedge^LX(t))|\\
			&\hspace{1em}=\binom{d-1}{L-1}\liminf_{t\to\infty}\frac{1}{t}\log|\det X(t)|\\
			&\hspace{1em}=\binom{d-1}{L-1}\sum_{i=1}^d\Lambda_i.
		\end{align*}
		Each index $i$ appears in precisely $\binom{d-1}{L-1}$ combinations of indices $i_1<\dots<i_L$. Hence, the above inequalities are actually equalities, proving that the induced basis is a normal basis and that the induced system is regular.
	\end{proof}
	
	The beginning of the proof of \cref{lemma:InducedNormalBasis} implies the following:
	
	\begin{proposition}\label[proposition]{prop:LE1Induced}
		It holds
		\begin{equation*}
			\lambda_{1,L}\leq \Lambda_1+\dots+\Lambda_L
		\end{equation*}
		independent of the regularity of \cref{eqn:system}.
	\end{proposition}
	
	\begin{proof}[proof of \cref{thm:METExtPowers}]
		Let $x_1,\dots,x_d$ be a normal basis such that $\chi[x_i]=\Lambda_i$. Since the induced basis is normal, it realizes the whole Lyapunov spectrum of the induced system. Our claims about the spectrum in \cref{thm:METExtPowers} now follow from
		$$\chi[x_{i_1}\wedge\dots\wedge x_{i_L}]=\Lambda_{i_1}\wedge\dots\wedge\Lambda_{i_L}.$$ 
		
		Moreover, $V_{2,L}$ is spanned by solutions
		\begin{equation*}
			x_{i_1}\wedge\dots\wedge x_{i_L}
		\end{equation*}
		such that 
		$$\Lambda_{i_1}+\dots+\Lambda_{i_L}<\Lambda_1+\dots+\Lambda_L.$$
		If $L=d_1+\dots+d_l$, the latter is equivalent to $x_{i_L}\in V_{l+1}$, which proves that $V_{2,L}$ is a subset of the set defined in \cref{thm:METExtPowers}. On the other hand, their dimensions must coincide since the codimension of $V_{2,L}$ is $d_{1,L}=1$ and neither set contains the solution $x_1\wedge\dots\wedge x_L$.
	\end{proof}
	
	As a direct consequence, we get \cref{prop:SingVals}.\par 
	
	\begin{proof}[proof of \cref{prop:SingVals}]
		The previous lemma and \cref{prop:CharExpMatrixNorm} imply
		\begin{align*}
			&\Lambda_L=(\Lambda_1+\dots+\Lambda_L)-(\Lambda_1+\dots+\Lambda_{L-1})\\
			&\hspace{1em}=\lim_{t\to\infty}\frac{1}{t}\log\|\wedge^LX(t)\|\\
			&\hspace{2em}-\lim_{t\to\infty}\frac{1}{t}\log\|\wedge^{L-1}X(t)\|\\
			&\hspace{1em}=\lim_{t\to\infty}\frac{1}{t}\log\frac{\sigma_1(X(t))\dots\sigma_L(X(t))}{\sigma_1(X(t))\dots\sigma_{L-1}(X(t))}\\
			&\hspace{1em}=\lim_{t\to\infty}\frac{1}{t}\log\sigma_L(X(t)).
		\end{align*}
	\end{proof}
	
	\begin{corollary}\label[corollary]{cor:StrongestGrowthInducedSystem}
		If \cref{eqn:system} is regular and $L=d_1+\dots+d_l$, then
		\begin{equation*}
			\lim_{t\to\infty} \frac{\|(\wedge^LX(t))\hat{v}\|}{\|\wedge^LX(t)\|}>0
		\end{equation*}
		for every $\hat{v}\notin V_{2,L,X}$.
	\end{corollary}
	
	\begin{proof}
		Since the induced system on the space of $L$-fold exterior products is regular and $\lambda_{1,L}$ is simple, the claim follows from \cref{cor:StrongestGrowth}.
	\end{proof}

\subsection{Stability of Lyapunov exponents}
	
	The stability of Lyapunov exponents has been widely studied. Here, we mainly state results listed in the overview from \cite[Ch.~V]{Adrianova1995IntroductionLinearSystems} and refer to \cite{BylovIzobov1969NecessarySufficientConditions} for the relation to integral separation.\par  
	
	\begin{definition}[{\cite[Def.~5.2.1]{Adrianova1995IntroductionLinearSystems}}]
		The Lyapunov exponents of \cref{eqn:system} are called \emph{stable} if for any $\epsilon>0$ there exist a $\delta>0$ such that the Lyapunov exponents of the continuously perturbed system $\dot{\tilde{x}}=(A(t)+Q(t))\tilde{x}$ with $\sup_{t\geq 0}\|Q(t)\|<\delta$ satisfy
		\begin{equation*}
			|\Lambda_i-\tilde{\Lambda}_i|<\epsilon
		\end{equation*}
		for $i=1,\dots,d$.
	\end{definition}
	
	The stability of Lyapunov exponents can be characterized by a uniform gap between the exponential growth rates of solutions corresponding to different Lyapunov exponents (this property is called \emph{integral separation}; see also \cite{BylovIzobov1969NecessarySufficientConditions}).
	
	\begin{theorem}[{\cite[Thm.~5.4.8]{Adrianova1995IntroductionLinearSystems}}]\label{thm:StabilityLE}
		The Lyapunov exponents of a system with simple spectrum are stable if and only if there exists a fundamental matrix $X(t)=[x_1(t),\dots,x_d(t)]$ and constants $a,b>0$ such that
		\begin{equation*}
			\frac{\|x_i(t)\|}{\|x_i(s)\|}\geq b e^{a(t-s)}\frac{\|x_{i+1}(t)\|}{\|x_{i+1}(s)\|}
		\end{equation*}
		for all $t\geq s$.
	\end{theorem}
	
	In fact, one may use such a fundamental matrix to construct a Lyapunov transformation that reduces \cref{eqn:system} to a diagonal system \cite[Thm.~5.3.1]{Adrianova1995IntroductionLinearSystems}.\par 
	
	The stability of Lyapunov exponents for systems with degenerate spectra can be characterized in a similar fashion by transforming the system to block diagonal form.
	
	\begin{theorem}[{\cite[Thm.~5.4.9]{Adrianova1995IntroductionLinearSystems},\cite{BarabanovDenisenko2007NecessarySufficientConditions}}]\label{thm:StabilityLE2}
		The Lyapunov exponents of \cref{eqn:system} are stable if and only if there exists a Lyapunov transformation reducing the system to block diagonal form
		\begin{equation*}
			\dot{y}=\textnormal{diag}(B_1(t),\dots,B_p(t))y,
		\end{equation*}
		where $B_i(t)\in\mathbb{R}^{d_i\times d_i}$ is upper triangular, such that the following hold:
		\begin{enumerate}[label=(\roman*)]
			\item all non-trivial solutions of $\dot{y}_i=B_i(t) y_i$ have characteristic exponent $\lambda_i$,
			\item $\lambda_i$ is stable\footnote{\cite[Thm.~5.4.9]{Adrianova1995IntroductionLinearSystems} uses the condition $\overline{\omega_i}=\lambda_i=\Omega_i$, which is equivalent to the stability of $\lambda_i$ on the subsystem (this follows from \cite[Thm.~5.4.9]{Adrianova1995IntroductionLinearSystems} applied to the subsystem).} for $\dot{y}_i=B_i(t) y_i$,
			\item there are constants $a,b>0$ such that
			\begin{equation}\label{eqn:UniformExpGap}
				\|Y_i(t,s)^{-1}\|^{-1}\geq b e^{a(t-s)}\|Y_{i+1}(t,s)\|
			\end{equation}
			for all $t\geq s$, where $Y_i(t,s)$ denotes the Cauchy matrix of $\dot{y}_i=B_i(t) y_i$.
		\end{enumerate}
	\end{theorem}	
	
	Systems with stable Lyapunov exponents retain their spectra under perturbations that tend to zero.\par 
	
	\begin{theorem}[{\cite{BarabanovDenisenko2007NecessarySufficientConditions}}]\label{thm:PerturbationToZero}
		Assume \cref{eqn:system} has stable Lyapunov exponents. If $Q(t)$ is a bounded, piecewise continuous perturbation such that
		\begin{equation*}
			\|Q(t)\|\to 0\text{ as }t\to\infty,
		\end{equation*}
		then the perturbed system $\dot{\tilde{x}}=(A(t)+Q(t))\tilde{x}$ has the same Lyapunov spectrum.
	\end{theorem}
	
	While $L^{\infty}$-perturbations of regular systems do not retain regularity in general, we have the following result:\par
	
	\begin{proposition}\label[proposition]{prop:PerturbationToZeroRegular}
		If \cref{eqn:system} is regular, has stable Lyapunov exponents, and $Q(t)$ is a bounded, piecewise continuous perturbation such that 
		\begin{equation*}
			\|Q(t)\|\to 0\text{ as }t\to\infty,
		\end{equation*}
		then the perturbed system is regular.
	\end{proposition}
	
	\begin{proof}
		\cref{thm:PerturbationToZero} implies that
		$$\Lambda_1+\dots+\Lambda_d=\tilde{\Lambda}_1+\dots+\tilde{\Lambda}_d.$$
		Moreover, due to the Liouville-Ostrogradski formula, we have
		\begin{align*}
			&\det \tilde{X}(t) \\
			&\hspace{1em}= \det\tilde{X}(0)\,e^{\int_{0}^{t}\textnormal{tr}(A(\tau)+Q(\tau))d\tau}\\
			&\hspace{1em}=\det (\tilde{X}(0) X(0)^{-1})\det X(t)\,e^{\int_{0}^{t}\textnormal{tr}(Q(\tau))d\tau}.
		\end{align*}
		Since 
		\begin{equation*}
			\lim_{t\to\infty}\frac{1}{t}\int_{0}^{t}\textnormal{tr}(Q(\tau))d\tau=0,
		\end{equation*}
		regularity of the original system implies regularity of the perturbed system:
		\begin{align*}
			\Lambda_1+\dots+\Lambda_d&=\tilde{\Lambda}_1+\dots+\tilde{\Lambda}_d\\
			&\geq \liminf_{t\to\infty}\frac{1}{t}\log|\det \tilde{X}(t)|\\
			&=\liminf_{t\to\infty}\frac{1}{t}\log|\det X(t)|\\
			&=\Lambda_1+\dots+\Lambda_d.
		\end{align*}
	\end{proof}

\subsection{Strong fast invertibility}
	
	In this subsection, we introduce \emph{strong fast invertibility}, which is a weaker concept than the stability of Lyapunov exponents but still sufficient for their computation. Our main objectives are a characterization result that allows us to compare strong fast invertibility to the stability of Lyapunov exponents and perturbation results in preparation for the analysis of Benettin's algorithm.\par

\subsubsection{Definition and relation to induced systems}

	Quas et al.\ introduce three notions of fast invertibility (weak/standard/strong) to analyze the existence of a dominated splitting for discrete-time systems \cite{QuasEtAl2019ExplicitBoundsSeparation}. Here, we focus solely on the strong version, since weak and standard fast invertibility are trivially satisfied for systems with bounded coefficients.\par 
	
	\begin{definition}[{\cite{QuasEtAl2019ExplicitBoundsSeparation}}]
		\cref{eqn:system} is said to be \emph{$L$-dimensionally strongly fast invertible}\footnote{The name ``fast invertibility'' stems from the context of \cite{QuasEtAl2019ExplicitBoundsSeparation} as it ensures that the cocylce is uniformly invertible on the fastest subspace of the respective dimension (without assuming bounded coefficients).} if
		\begin{equation*}
			c_{FI,L}:=\sup_{t\geq s\geq \tau}\prod_{k=1}^{L}\frac{\sigma_k(X(t,s))\sigma_k(X(s,\tau))}{\sigma_k(X(t,\tau))}<\infty.
		\end{equation*}
	\end{definition}
	
	There are several equivalent formulations immediately visible from the definition.\par 
	
	\begin{lemma}\label[lemma]{lemma:FIExtPowers}
		The following are equivalent:
		\begin{enumerate}[label=(\roman*)]
			\item \cref{eqn:system} is $L$-dim.\ strongly fast invertible.
			\item The induced system on $\wedge^L\mathbb{R}^d$ is $1$-dim.\ strongly fast invertible.
			\item There is a constant $c>0$ such that
			\fontsize{8pt}{0pt}
			\begin{equation*}
				\frac{\|\wedge^LX(t,\tau)\|}{\|\wedge^LX(s,\tau)\|}\leq \|\wedge^L X(t,s)\|\leq c\frac{\|\wedge^LX(t,\tau)\|}{\|\wedge^LX(s,\tau)\|}
			\end{equation*}
			\normalsize
			for all $t\geq s\geq \tau$.
			\item $\sup_{t\geq s}\prod_{k=1}^{L}\frac{\sigma_k(X(t,s))\sigma_k(X(s))}{\sigma_k(X(t))}<\infty$ for any fundamental matrix $X(t)$.
			\item $\sup_{t\geq s\geq \tau\geq T}\prod_{k=1}^{L}\frac{\sigma_k(X(t,s))\sigma_k(X(s,\tau))}{\sigma_k(X(t,\tau))}<\infty$ for any $T>0$.
		\end{enumerate}
	\end{lemma}
	
	\begin{proof}
		The relation between strong fast invertibility and exterior products follows from \cref{lemma:ExteriorPowerMatrixProperties}. The third characterization is merely a reformulation of strong fast invertibility for the induced system. Equivalence of strong fast invertibility to the fourth characterization follows from
		\begin{align*}
			&\frac{\|\wedge^LX(t,s)\|\,\|\wedge^LX(s)\|}{\|\wedge^LX(t)\|}\\
			&\hspace{1em}\leq \|\wedge^LX(0)\|\,\|\wedge^LX(0)^{-1}\|\\
			&\hspace{2em}\cdot\frac{\|\wedge^LX(t,s)\|\,\|\wedge^LX(s,0)\|}{\|\wedge^LX(t,0)\|}
		\end{align*}
		and
		\begin{align*}
			&\frac{\|\wedge^LX(t,s)\|\,\|\wedge^LX(s,\tau)\|}{\|\wedge^LX(t,\tau)\|}\\
			&\hspace{1em}= \frac{\|\wedge^LX(t,s)\|\,\|\wedge^LX(s)\|}{\|\wedge^LX(t)\|}\\
			&\hspace{2em}\cdot\frac{\|\wedge^LX(s,\tau)\|\,\|\wedge^LX(\tau)\|}{\|\wedge^LX(s)\|}\\
			&\hspace{2em}\cdot\frac{\|\wedge^LX(t)\|}{\|\wedge^LX(t,\tau)\|\,\|\wedge^LX(\tau)\|}\\
			&\hspace{1em}\leq \frac{\|\wedge^LX(t,s)\|\,\|\wedge^LX(s)\|}{\|\wedge^LX(t)\|}\\
			&\hspace{2em}\cdot\frac{\|\wedge^LX(s,\tau)\|\,\|\wedge^LX(\tau)\|}{\|\wedge^LX(s)\|}.
		\end{align*}
		Finally, using
		\begin{align*}
			&\wedge^L{X}(T+t,t)\wedge^L{X}(t,s)\\
			&\hspace{1em}=\wedge^L{X}(T+t,T+s)\wedge^L{X}(T+s,s),
		\end{align*}
		and $\sup_t\|\hat{A}^L(t)\|\leq LM$, one may show the final characterization:
		\small
		\begin{align*}
			&\frac{\|\wedge^L{X}(t,s)\|\,\|\wedge^L{X}(s,\tau)\|}{\|\wedge^L{X}(t,\tau)\|}\leq e^{6TLM}\\
			&\hspace{2em}\cdot\frac{\|\wedge^L{X}(T+t,T+s)\|\,\|\wedge^L{X}(T+s,T+\tau)\|}{\|\wedge^L{X}(T+t,T+\tau)\|}.
		\end{align*}
		\normalsize	
	\end{proof}
	
	Maybe the most intuitive of the above characterizations of strong fast invertibility is the third. It says that the maximal growth of $L$-dim.\ volumes is (up to a constant) the same as the maximal growth of $L$-dim.\ volumes relative to any initial time.\par 
	
	\begin{proposition}
		Every system is $d$-dim.\ strongly fast invertible.
	\end{proposition}
	
	\begin{proof}
		Since $\dim(\wedge^d\mathbb{R}^d)=1$, we have $c_{FI,d}=1$.
	\end{proof}

\subsubsection{Characterization and comparison to stability of Lyapunov exponents}
	
	To derive the characterization theorem mentioned in the introduction, we first prove that regularity and strong fast invertibility imply the existence of a Lyapunov transformation that brings the system into block diagonal form.\par 
	
	\begin{proposition}\label[proposition]{prop:FIAngleGap}
		Assume \cref{eqn:system} has Lyapunov exponents $\lambda_1>\dots>\lambda_p$ with multiplicities $d_1+\dots+d_p=d$. If the system is regular and strongly fast invertible at dim.\ $d_1+\dots+d_l$ for $l=1,\dots,p$, then there is a Lyapunov transformation reducing the system to block diagonal form:
		\begin{equation*}
			\dot{y}=\textnormal{diag}(B_1(t),\dots,B_p(t))y,
		\end{equation*}
		where $B_i\in\mathbb{R}^{d_i\times d_i}$ is upper triangular, such that all non-trivial solutions of $\dot{y}_i=B_i y_i$ have characteristic exponent $\lambda_i$.
	\end{proposition}
	
	We prove the proposition via induced systems using two auxiliary lemmata and the following result:

	\begin{theorem}[{\cite[Thm.~3.3.3 and Thm.~3.3.4]{Adrianova1995IntroductionLinearSystems}}]\label[theorem]{thm:BlockTriang}
		There exists a Lyapunov transformation reducing \cref{eqn:system} to block diagonal form:
		\begin{equation*}
			\dot{y}=\textnormal{diag}(B_1(t),\dots,B_k(t))y,
		\end{equation*}
		where $B_i\in\mathbb{R}^{n_i\times n_i}$ is upper triangular, if and only if there is a fundamental matrix 
		$$X=[X_{1},\dots,X_{k}]$$
		with $X_i\in\mathbb{R}^{d\times n_i}$ satisfying
		\begin{equation}\label{eqn:AngleGap}
			\inf_t\frac{G(X)}{G(X_{1})\dots G(X_{k})}>0,
		\end{equation}
		where $G$ denotes the Gram determinant.
	\end{theorem}

	\begin{remark}\label[remark]{remark:BlockTriang}
		Given a fundamental matrix $X(t)$ as in \cref{thm:BlockTriang}, the Lyapunov transformation is constructed by applying the Gram-Schmidt procedure to each block $X_{i}(t)$ individually, yielding 
		\begin{align*}
			X(t)&=[Q_{1}(t),\dots,Q_{k}(t)]\textnormal{diag}(R_{1}(t),\dots,R_{k}(t))\\
			&:=L(t)Y(t).
		\end{align*} 
	\end{remark}
	
	\begin{lemma}[{\cite[Remark~3.3.4]{Adrianova1995IntroductionLinearSystems}}]\label[lemma]{lemma:FIAnglePrincipleAngles}
		Let $X=[X_1,X_2]$. It holds
		\begin{equation*}
			\frac{G(X)}{G(X_1)G(X_2)}=\prod_{i=1}^{l}\sin^2 \alpha_i(\textnormal{im}X_1,\textnormal{im}X_2),
		\end{equation*}
		where $\alpha_1\leq \dots\leq \alpha_l$ are the principle angles.
	\end{lemma}
	
	\begin{lemma}\label[lemma]{lemma:FIAngleGapReducedForm}
		If \cref{eqn:system} is $1$-dim.\ strongly fast invertible and $X=[X_1,X_2]$ is a fundamental matrix such that $X_1\in\mathbb{R}^{d\times 1}$ has a higher characteristic exponent than the columns of $X_2\in\mathbb{R}^{d\times (d-1)}$, then 
		\begin{equation}\label{eqn:Alpha1}
			\inf_t\sin^2\alpha_1(\textnormal{im}X_1(t),\textnormal{im}X_2(t))>0.
		\end{equation}
	\end{lemma}
	
	\begin{proof}
		Let $X=[X_1,X_2]$ be as in the claim. We first apply the Gram-Schmidt procedure to bring $X$ into upper triangular form. Indeed, the Gram-Schmidt procedure is an orthogonal Lyapunov transformation (\cref{proposition:GramSchmidt}).\ Thus, it leaves the fraction in \cref{eqn:AngleGap} and hence also $\sin^2\alpha_1(\textnormal{im}X_1,\textnormal{im}X_2)$ invariant. So, it is sufficient to check \cref{eqn:Alpha1} for
		\begin{equation*}
			Y=[Y_1,Y_2]=\begin{pmatrix}
				y_{11} & Y_{12}\\
				0 & Y_{22}
			\end{pmatrix}
		\end{equation*}
		in which $y_{11}$ has a higher characteristic exponent than $\|Y_{12}\|$ and $\|Y_{22}\|$.\par 
		
		Given $s\geq 0$, we always find $t_s\geq s$ such that
		\begin{equation*}
			\max_{\alpha\neq 0}\frac{\left|\frac{y_{11}(s)}{y_{11}(t_s)}Y_{12}(t_s)\alpha\right|}{\|Y_{22}(s)\alpha\|}\leq 1
		\end{equation*}
		and
		\begin{equation*}
			\|Y(t_s)\|\leq 2|y_{11}(t_s)|.
		\end{equation*}
		In order to make use of strong fast invertibility, we compute the Cauchy matrix:
		\small
		\begin{align*}
			&Y(t,s)=\begin{pmatrix}
				y_{11}(t,s) & Y_{12}(t,s)\\
				0 & Y_{22}(t,s)
			\end{pmatrix}\\
			&=\begin{pmatrix}
				\frac{y_{11}(t)}{y_{11}(s)} & Y_{12}(t)Y_{22}(s)^{-1}-\frac{y_{11}(t)}{y_{11}(s)}Y_{12}(s)Y_{22}(s)^{-1}\\
				0 & Y_{22}(t)Y_{22}(s)^{-1}
			\end{pmatrix}.
		\end{align*}
		\normalsize
		Now, for $s\geq 0$ and $\alpha\neq 0$ it holds
		\begin{align*}
			&\frac{|Y_{12}(s)\alpha|}{\|Y_{22}(s)\alpha\|}\\
			&\hspace{1em}\leq 1 + \left|\frac{\frac{y_{11}(s)}{y_{11}(t_s)}Y_{12}(t_s)\alpha}{\|Y_{22}(s)\alpha\|} - \frac{ Y_{12}(s)\alpha}{\|Y_{22}(s)\alpha\|}\right|\\
			&\hspace{1em}=1+\frac{\left|\frac{y_{11}(s)}{y_{11}(t_s)}Y_{12}(t_s,s)Y_{22}(s)\alpha\right|}{\|Y_{22}(s)\alpha\|}\\
			&\hspace{1em}\leq 1+\frac{\|Y_{12}(t_s,s)\|\,|y_{11}(s)|}{|y_{11}(t_s)|}\\
			&\hspace{1em}\leq 1+2\frac{\|Y(t_s,s)\|\,\|Y(s)\|}{\|Y(t_s)\|}\\
			&\hspace{1em}\leq \underset{=:c}{\underbrace{1+2c_{FI,L}\,\|Y(0)\|\,\|Y(0)^{-1}\|}}.\\
		\end{align*}
		Since
		\begin{align*}
			\left(\frac{\left \langle Y_1(t),Y_2(t)\alpha \right\rangle}{\left\|Y_1(t)\right\|\,\left\|Y_2(t)\alpha\right\|}\right)^2&=\left(\frac{\left \langle e_1,Y_2(t)\alpha \right\rangle}{\left\|Y_2(t)\alpha\right\|}\right)^2\\
			&=\frac{|Y_{12}(t)\alpha|^2}{|Y_{12}(t)\alpha|^2+\|Y_{22}(t)\alpha\|^2}\\
			&=\frac{1}{1+\left(\frac{|Y_{12}(t)\alpha|}{\|Y_{22}(t)\alpha\|}\right)^{-2}}\\
			&\leq \frac{1}{1+c^{-2}},
		\end{align*}
		we have
		\begin{align*}
			&\sin^2\alpha_1(\textnormal{im}Y_1,\textnormal{im}Y_2)\\
			&\hspace{1em}=1-\cos^2\alpha_1(\textnormal{im}Y_1,\textnormal{im}Y_2)\\
			&\hspace{1em}=1-\max_{\alpha\neq 0}\left(\frac{\left \langle Y_1(t),Y_2(t)\alpha \right\rangle}{\left\|Y_1(t)\right\|\,\left\|Y_2(t)\alpha\right\|}\right)^2\\
			&\hspace{1em}\geq 1-\frac{1}{1+c^{-2}}\\
			&\hspace{1em}>0.
		\end{align*}
	\end{proof}
	
	\begin{proof}[proof of \cref{prop:FIAngleGap}]
		Let $X(t)=[X_{1},\dots,X_{p}]$ be an ordered normal fundamental matrix so that the columns of $X_{i}$ have characteristic exponent $\lambda_i$. We first show that \cref{eqn:AngleGap} holds with $n_1=L$ and $n_2=d-L$ if the system is strongly fast invertible at dim.\ $L=d_1+\dots+d_l$ and then apply \cref{thm:BlockTriang} successively to arrive at the desired form.\par
		
		Fix $L=d_1+\dots+d_l$ and set
		$$U:=\textnormal{im}[X_1,\dots,X_l]\text{ and }V:=\textnormal{im}[X_{l+1},\dots,X_p].$$	
		\cref{lemma:FIAnglePrincipleAngles} and \cref{cor:PrincipleAngles} imply
		\begin{align*}
			&\frac{G(X)}{G([X_{1},\dots,X_{l}])G([X_{l+1},\dots,X_{p}])}\\
			&\hspace{1em}=\prod_{i=1}^{\min(L,d-L)}\sin^2\alpha_i(U,V)\\
			&\hspace{1em}=\sin^2\alpha_1(\hat{U},\hat{V}),
		\end{align*}		
		where
		\begin{equation*}
			\hat{U}=\textnormal{span}\{(\wedge^LX)(e_1\wedge\dots\wedge e_L)\}
		\end{equation*}
		and
		\begin{align*}
			\hat{V}&=\textnormal{span}\{(\wedge^LX)(e_{i_1}\wedge\dots\wedge e_{i_L})\ |\\
			&\hspace{2em}\{i_1,\dots,i_L\}\neq \{1,\dots,L\}\}.
		\end{align*}
		Now, let $Y=[Y_1,Y_2]$ be a matrix representation of $\wedge^LX$ with respect to the standard basis on $\wedge^L\mathbb{R}^d$. Order the induced basis starting with $e_1\wedge\dots \wedge e_L$ so that the first column $Y_1$ of $Y$ represents $(\wedge^LX)(e_1\wedge\dots\wedge e_L)$. Then, 
		\begin{equation*}
			\sin^2\alpha_1(\hat{U},\hat{V})=\sin^2\alpha_1(\textnormal{im}Y_1,\textnormal{im}Y_2).
		\end{equation*}
		Due to regularity, $Y_1$ has characteristic exponent $\Lambda_1+\dots+\Lambda_L$ and the columns of $Y_2$ have lower characteristic exponents. Hence, \cref{lemma:FIAngleGapReducedForm} and the above arguments imply
		\begin{equation*}
			\inf_t\frac{G(X)}{G([X_{1},\dots,X_{l}])G([X_{l+1},\dots,X_{p}])}>0.
		\end{equation*}
				
		To arrive at the block diagonal form claimed in the proposition, we apply the transformation from \cref{remark:BlockTriang} successively. First, we use strong fast invertibility at dim.\ $d_1$ to make the first $d_1$ columns of $X$ orthogonal to the rest. Assuming orthogonality of the first $d_1$ columns, strong fast invertibility at dim.\ $d_1+d_2$ implies
		\begin{align*}
			0&<\inf_t\frac{G(X)}{G([X_{d_1},X_{d_2}])G([X_{d_3},\dots,X_{d_p}])}\\
			&=\inf_t\frac{G(X)}{G(X_{d_1})G(X_{d_2})G([X_{d_3},\dots,X_{d_p}])}.
		\end{align*}
		Thus, we may apply another Lyapunov transformation to decouple the next subsystem and so on.
	\end{proof}
	
	Using \cref{prop:FIAngleGap}, we now may prove our characterization result for strong fast invertibility.\par 
	
	\begin{theorem}\label[theorem]{thm:FICharacterization}
		Assume \cref{eqn:system} has Lyapunov exponents $\lambda_1>\dots>\lambda_p$ with multiplicities $d_1+\dots+d_p=d$. If the system is regular and strongly fast invertible at dim.\ $d_1+\dots+d_l$ for $l=1,\dots,p$, then there exists a Lyapunov transformation reducing the system to block diagonal form
		\begin{equation*}
			\dot{y}=\textnormal{diag}(B_1(t),\dots,B_p(t))y,
		\end{equation*}
		where $B_i(t)\in\mathbb{R}^{d_i\times d_i}$ is upper triangular, such that the following hold:
		\begin{enumerate}[label=(\roman*)]
			\item all non-trivial solutions of $\dot{y}_i=B_i(t) y_i$ have characteristic exponent $\lambda_i$,
			\item there is a constant $b>0$ such that
				\begin{equation}\label{eqn:UniformGap}
					\|Y_i(t,s)^{-1}\|^{-1}\geq b\|Y_{i+1}(t,s)\|
				\end{equation}
				for all $t\geq s$, where $Y_i(t,s)$ denotes the Cauchy matrix of $\dot{y}_i=B_i(t) y_i$.
		\end{enumerate}
		Conversely, any block diagonal system $\dot{y}=B(t)y$ satisfying $(i)$ and $(ii)$ is strongly fast invertible at dim.\ $d_1+\dots+d_l$ for $l=1,\dots,p$.
	\end{theorem}  

	\begin{proof}
		Let us first assume that the system is regular and strongly fast invertible at the respective dimensions. After applying \cref{prop:FIAngleGap}, all that remains is to show that \cref{eqn:UniformGap} is satisfied for the reduced system. To this end, fix $L=d_1+\dots+d_l$ and set $Y=Z_1+Z_2$ with
		\begin{equation*}
			Z_1 := \textnormal{diag}(Y_1,\dots,Y_L,0,\dots,0)
		\end{equation*}
		and
		\begin{equation*}
			Z_2 := \textnormal{diag}(0,\dots,0,Y_{L+1},\dots,Y_p).
		\end{equation*}
		Due to \cref{lemma:ExtPowerPropertiesDiag}, we have
		\begin{equation*}
			\|\wedge^LY\| = \max_{k}\|(\wedge^kZ_1)\wedge(\wedge^{L-k}Z_2)\|
		\end{equation*}
		and
		\begin{align*}
			\|\wedge^LY\|\geq \|\wedge^{L-1}Z_1\|\,\|Z_2\|=\frac{\|\wedge^LZ_1\|}{\sigma_L(Z_1)}\|Z_2\|.
		\end{align*}
		The singular value $\sigma_L(Z_1)$ can be estimated via
		\begin{align*}
			\sigma_L(Z_1)&=\sigma_L(\textnormal{diag}(Y_1,\dots,Y_l))\\
			&=\|\textnormal{diag}(Y_1^{-1},\dots,Y_l^{-1})\|^{-1}\\
			&\leq \|Y_l^{-1}\|^{-1}.
		\end{align*}
		Due to regularity, we have 
		\small
		$$\chi[\|\wedge^LZ_1\|]=\lambda_{1,L}>\lambda_{2,L}\geq \chi[\|(\wedge^kZ_1)\wedge(\wedge^{L-k}Z_2)\|]$$
		\normalsize
		for $k<L$, which implies
		\begin{equation*}
			\|\wedge^LY(t)\|=\|\wedge^LZ_1(t)\|
		\end{equation*}
		for $t$ large enough. Moreover, it holds
		\begin{align*}
			\|\wedge^LZ_1(t,s)\|&=\|\wedge^L\textnormal{diag}(Y_1(t,s),\dots,Y_l(t,s))\|\\
			&=\frac{\|\wedge^L\textnormal{diag}(Y_1(t),\dots,Y_l(t))\|}{\|\wedge^L\textnormal{diag}(Y_1(s),\dots,Y_l(s))\|}\\
			&=\frac{\|\wedge^LZ_1(t)\|}{\|\wedge^LZ_1(s)\|}
		\end{align*}
		since $\wedge^L\mathbb{R}^L$ is one-dimensional.\par
		
		Combining the previous estimates, we get
		\begin{align*}
			&\|Y_l(t,s)^{-1}\|^{-1}\\
			&\hspace{1em}\geq\frac{\|\wedge^L Z_1(t,s)\|}{\|\wedge^LY(t,s)\|}\|Z_2(t,s)\|\\
			&\hspace{1em}\geq  \frac{\|\wedge^L Z_1(t,s)\|}{\|\wedge^LY(t,s)\|}\|Y_{l+1}(t,s)\|\\
			&\hspace{1em}=  \frac{\|\wedge^L Z_1(t)\|}{\|\wedge^LY(t,s)\|\,\|\wedge^LZ_1(s)\|}\|Y_{l+1}(t,s)\|\\
			&\hspace{1em}\geq  \frac{\|\wedge^L Y(t)\|}{\|\wedge^LY(t,s)\|\,\|\wedge^LY(s)\|}\|Y_{l+1}(t,s)\|\\
			&\hspace{1em}\geq  \frac{1}{c_{FI,L}\|\wedge^LY(0)\|\,\|\wedge^LY(0)^{-1}\|}\|Y_{l+1}(t,s)\|
		\end{align*}
		for large $t>0$, which proves \cref{eqn:UniformGap}.\par 
		
		Next, we show that any block diagonal system as in \cref{thm:FICharacterization} is strongly fast invertible at the respective dimensions. As before, fix $L=d_1+\dots+d_l$ and assume the decomposition $Y=Z_1+Z_2$. Applying \cref{eqn:UniformGap} repeatedly yields
		\begin{align*}
			\|Y_i^{-1}\|^{-1}&\geq b\|Y_{i+1}\|\geq b\|Y_{i+1}^{-1}\|^{-1}\geq\dots\\
			&\geq b^{k}\|Y_{i+k}\|
		\end{align*}
		and hence	
		\begin{align*}
			\sigma_L(Z_1)&=\sigma_L(\textnormal{diag}(Y_1,\dots,Y_l))\\
			&=\|\textnormal{diag}(Y_1^{-1},\dots,Y_l^{-1})\|^{-1}\\
			&= \min_{i=1,\dots,l}\|Y_i^{-1}\|^{-1}\\
			&\geq b^{p-1}\max_{i=l+1,\dots,p}\|Y_i\|\\
			&= b^{p-1}\|Z_2\|
		\end{align*}		
		since $b\in(0,1]$ (set $s=t$ in \cref{eqn:UniformGap}). Now, we estimate
		\begin{align*}
		&\|\wedge^LY\|=\max_k\|(\wedge^kZ_1)\wedge(\wedge^{L-k}Z_2)\|\\
		&\hspace{1em}\leq \max_k\binom{d}{L}^{\frac{1}{2}}\|\wedge^kZ_1\|\,\|\wedge^{L-k}Z_2\|\\
		&\hspace{1em}\leq \max_k\binom{d}{L}^{\frac{1}{2}}\|\wedge^kZ_1\|\,\|Z_2\|^{L-k}\\
		&\hspace{1em}\leq \max_k\binom{d}{L}^{\frac{1}{2}}\|\wedge^kZ_1\|b^{-(p-1)(L-k)}\sigma_L(Z_1)^{L-k}\\
		&\hspace{1em}\leq \underset{=:c}{\underbrace{\max_k\left(\binom{d}{L}^{\frac{1}{2}}b^{-(p-1)(L-k)}\right)}}\|\wedge^LZ_1\|.
		\end{align*}
		Hence, we get
		\begin{align*}
			&\sup_{t\geq s}\frac{\|\wedge^LY(t,s)\|\,\|\wedge^LY(s)\|}{\|\wedge^LY(t)\|}\\
			&\hspace{1em}\leq c^2\underset{=1}{\underbrace{\sup_{t\geq s}\frac{\|\wedge^L Z_1(t,s)\|\,\|\wedge^LZ_1(s)\|}{\|\wedge^LZ_1(t)\|}}}<\infty.
		\end{align*}
	\end{proof}
	
	\begin{remark}
		When testing for strong fast invertibility (resp.\ for stability of Lyapunov exponents), it is enough to check the conditions of \cref{thm:FICharacterization} (resp.\ \cref{thm:StabilityLE2}) after transforming the original system to any block diagonal system of the form
		\begin{equation*}
			\dot{y}=\textnormal{diag}(B_1(t),\dots,B_p(t))y
		\end{equation*}
		such that the non-trivial solutions of block $B_i(t)\in\mathbb{R}^{d_i\times d_i}$ have characteristic exponent $\lambda_i$.
	\end{remark}
	
	\begin{proof}
		 Any Lyapunov transformation between systems of this form acts as Lyapunov transformations on the individual blocks, i.e., $L(t)=\textnormal{diag}(L_1(t),\dots,L_p(t))$. This implies that \cref{eqn:UniformGap} (resp.\ stability of the Lyapunov exponents of the subsystems and \cref{eqn:UniformExpGap}) holds for one such system if and only if it holds for all of them.
	\end{proof}
	
	An immediate consequence of \cref{thm:FICharacterization} is that stability of Lyapunov exponents implies strong fast invertibility.\par 
	
	\begin{corollary}
		Assume \cref{eqn:system} has Lyapunov exponents $\lambda_1>\dots>\lambda_p$ with multiplicities $d_1+\dots+d_p=d$. If the Lyapunov exponents are stable, then the system is strongly fast invertible at dim.\ $d_1+\dots+d_l$ for $l=1,\dots,p$.
	\end{corollary}
	
	\begin{remark}
		There is a regular system that is strongly fast invertible at every dimension but does not have stable Lyapunov exponents.
	\end{remark}
	
	\begin{proof}
	Let $\epsilon>0$. We set $a_n:=n^2-\log n$ and $b_n:=n^2$ and define $f_n\in C(\mathbb{R}_{\geq 0},\mathbb{R})$ through
	\begin{equation*}
		f_n=\begin{cases}
			0,\quad &t\in [0,a_n-\epsilon]\\
			\frac{t-(a_n-\epsilon)}{\epsilon}, &t\in (a_n-\epsilon,a_n]\\
			1,\quad &t\in (a_n,b_n]\\
			\frac{b_n+\epsilon-t}{\epsilon}, &t\in (b_n,b_n+\epsilon]\\
			0. &t\in (b_n+\epsilon,\infty)
		\end{cases}
	\end{equation*}
	For small $\epsilon$ the functions $f_n$ have disjoint supports. Thus, the system matrix
	\begin{equation*}
		A(t):=\begin{pmatrix}
			1 & 0\\
			0 & a_{22}(t)
		\end{pmatrix}
	\end{equation*}
	with $a_{22}(t):=\sum_{n\in\mathbb{N}}f_n(t)$ is bounded and continuous. A fundamental matrix is given by
	\begin{equation*}
		X(t):=\begin{pmatrix}
			e^t & 0\\
			0 & e^{\int_0^t a_{22}(\tau)\,d\tau}
		\end{pmatrix}.
	\end{equation*}
	Since
	\begin{equation*}
		1\leq e^{\int_{0}^{b_n+\epsilon}a_{22}(\tau)\,d\tau}\leq e^{\sum_{k=1}^{n}(\log(k) + 2\epsilon)}=n!e^{2\epsilon n},
	\end{equation*}
	we have
	\begin{equation*}
		0\leq \frac{1}{t}\log e^{\int_{0}^{t}a_{22}(\tau)\,d\tau}\leq \frac{2\epsilon (n+1)}{n^2} + \frac{\log((n+1)!)}{n^2} 
	\end{equation*}
	for $t\in [b_n,b_{n+1}+\epsilon]$. In particular, the system has Lyapunov exponents $\lambda_1=1$ and $\lambda_2=0$ and is regular. Moreover, since 
	\begin{equation*}
		\|X(t,s)\|_{\max}=e^{t-s},
	\end{equation*}
	the system is strongly fast invertible at every dimension.\par 
	
	However, the Lyapunov exponents are not stable, because there are no constants $a,b>0$ such that
	\begin{equation*}
		e^{(b_n-a_n)-\int_{a_n}^{b_n} a_{22}(\tau)d\tau}=1
	\end{equation*}
	is bounded from below by $be^{a(b_n-a_n)}=bn^a$ for all $n\in\mathbb{N}$.
	\end{proof}
	
	In case the Lyapunov spectrum is simple, strong fast invertibility implies the existence of a fundamental system of solutions with uniformly separated growth.\par
	
	\begin{theorem}\label[theorem]{thm:FICharacterizationSimple}
			If \cref{eqn:system} is regular, has simple Lyapunov spectrum and is strongly fast invertible at every dimension, then there exists a fundamental matrix $X(t)=[x_1(t),\dots,x_d(t)]$ and a constant $b>0$ such that
			\begin{equation}\label{eqn:Separation}
				\frac{\|x_i(t)\|}{\|x_i(s)\|}\geq b\frac{\|x_{i+1}(t)\|}{\|x_{i+1}(s)\|}
			\end{equation}
			for all $t\geq s$.
	\end{theorem}
	
	\begin{proof}
		According to \cref{thm:FICharacterization} there is a Lyapunov transformation $L(t)$ bringing the system to diagonal form with fundamental matrix $Y=\textnormal{diag}(y_{11},\dots,y_{dd})$ and a constant $b'>0$ such that
		$$\frac{|y_{ii}(t)|}{|y_{ii}(s)|}\geq b'\frac{|y_{(i+1)(i+1)}(t)|}{|y_{(i+1)(i+1)}(s)|}$$ 
		for all $t\geq s$. The fundamental matrix $X=LY$ of the original system satisfies
		$$\|x_i(t)\|=\|L(t)e_i\|\,|y_{ii}(t)|.$$
		Since $L(t)$ is a Lyapunov transformation, the term $\|L(t)e_i\|$ is uniformly bounded from below and above by positive constants. \cref{eqn:Separation} follows.
	\end{proof}
	
	\begin{remark}
		There is a regular system with simple Lyapunov spectrum satisfying \cref{eqn:Separation} that is not strongly fast invertible at every dimension.
	\end{remark}
	
	\begin{proof}
		Let $0<\epsilon_n\leq 1/e$ be a sequence of numbers converging to zero. We set 
		$$c_n:=\left(\frac{9}{4}\right)^n\prod_{k=1}^{n}\left(\frac{1}{\epsilon_k}\right)^{\frac{1}{\epsilon_k}}.$$
		The numbers
		\begin{align*}
			t_{n,0}&:=c_{n},\\
			t_{n,1}&:=c_{n}+\frac{2}{3(1-\epsilon_n)},\\
			t_{n,2}&:=c_{n}+\frac{2}{3(1-\epsilon_n)}+\frac{\log\left(\frac{1}{\epsilon_n}\right)}{1-\epsilon_n}-\frac{2}{3},\\
			t_{n,3}&:=c_{n}+\frac{2}{3(1-\epsilon_n)}+\frac{\log\left(\frac{1}{\epsilon_n}\right)}{1-\epsilon_n},\\
			t_{n,4}&:=c_{n}+\frac{2}{3(1-\epsilon_n)}+\frac{\log\left(\frac{1}{\epsilon_n}\right)}{\epsilon_n(1-\epsilon_n)},\\
			t_{n,5}&:=c_{n}+\frac{2}{3(1-\epsilon_n)}+\frac{\log\left(\frac{1}{\epsilon_n}\right)}{\epsilon_n(1-\epsilon_n)}+1,\\
			t_{n,6}&:=c_{n}+\frac{5}{3(1-\epsilon_n)}+\frac{\log\left(\frac{1}{\epsilon_n}\right)}{\epsilon_n(1-\epsilon_n)},
		\end{align*}
		satisfy
		$$t_{n,0}\leq \dots\leq t_{n,6}\leq t_{n+1,0}.$$
		On the partition introduced by these numbers, we define functions
		\begin{equation*}
			f(t):=\begin{cases}
				1, &t\in[0,t_{1,0})\\
				c_{n-1}p_{n,f}(t-t_{n,0}), &t\in[t_{n,0},t_{n,1})\\
				c_{n-1}\frac{3}{2}e^{(1-\epsilon_n)(t-t_{n,1})}, &t\in[t_{n,1},t_{n,4})\\
				c_{n-1}\frac{3}{2}\left(\frac{1}{\epsilon_n}\right)^{\frac{1}{\epsilon_n}}\\
				\hspace{1em}\cdot q_{n,f}(t-t_{n,4}), &t\in[t_{n,4},t_{n,6})\\
				c_{n}, &t\in[t_{n,6},t_{n+1,0})
			\end{cases}
		\end{equation*}
		and
		\begin{equation*}
			g(t):=\begin{cases}
				1, &t\in[0,t_{1,2})\\
				c_{n-1}p_{n,g}(t-t_{n,2}), &t\in[t_{n,2},t_{n,3})\\
				c_{n-1}\frac{3}{2}e^{t-t_{n,3}}, &t\in[t_{n,3},t_{n,4})\\
				c_{n-1}\frac{3}{2}\left(\frac{1}{\epsilon_n}\right)^{\frac{1}{\epsilon_n}}\\
				\hspace{1em}\cdot q_{n,g}(t-t_{n,4}), &t\in[t_{n,4},t_{n,5})\\
				c_{n}, &t\in[t_{n,5},t_{n+1,2})
			\end{cases}
		\end{equation*}
		using the polynomials
		\begin{align*}
			p_{n,f}(t)&:=\frac{9}{8}(1-\epsilon_n)^2t^2+1,\\
			p_{n,g}(t)&:=\frac{9}{8}t^2+1,\\
			q_{n,f}(t)&:=-\frac{1}{2}(1-\epsilon_n)^2t^2+(1-\epsilon_n)t+1,\\
			q_{n,g}(t)&:=-\frac{1}{2}t^2+t+1,
		\end{align*}
		for smoothing. The functions have the following properties:
		\begin{itemize}
			\item $f,g\in C^1(\mathbb{R}_{\geq 0},\mathbb{R})$,
			\item $f,g$ are monotonically increasing,
			\item $|f'-f|/|g|$ bounded,
			\item $|g'|/|g|$ bounded,
			\item $1\leq f(t),g(t)\leq t$ for $t\geq 1$.
		\end{itemize}
		In particular,
		\begin{equation*}
			\dot{x}=A(t)x:=\begin{pmatrix}
				1 & \frac{f'-f}{g}\\
				0 & \frac{g'}{g}
			\end{pmatrix}x
		\end{equation*}
		is a continuous, bounded system with fundamental matrix 
		\begin{equation*}
			X(t)=[x_1(t),x_2(t)]=\begin{pmatrix}
				e^t & f \\
				0 & g
			\end{pmatrix}.
		\end{equation*} 
		The solution $x_2(t)$ has characteristic exponent zero and the system is regular.\par
		
		To show property \cref{eqn:Separation}, we first remark that
		\begin{equation*}
			\frac{|p_{n,f}'|}{|p_{n,f}|},\frac{|q_{n,f}'|}{|q_{n,f}|}\leq 1
		\end{equation*} 
		on the respective intervals used in the definition of $f$. Since 
		\begin{equation*}
			(\log f(t))'=\frac{f'(t)}{f(t)}\leq 1,
		\end{equation*}
		the mean value theorem implies
		\begin{equation*}
			\frac{f(t)}{f(s)}\leq e^{t-s}
		\end{equation*}
		for $t\geq s$. In combination with $g(t)\leq \frac{3}{2} f(t)$, we get \cref{eqn:Separation}:
		\begin{align*}
			\frac{\|x_2(t)\|^2}{\|x_2(s)\|^2}&=\frac{f^2(t)+g^2(t)}{f^2(s)+g^2(s)}\\
			&\leq \frac{13}{4}\frac{f^2(t)}{f^2(s)} \\
			&\leq \frac{13}{4}\frac{\|x_1(t)\|^2}{\|x_1(s)\|^2}
		\end{align*}
		for $t\geq s$.\par 
		
		To show that the system is not $1$-dim.\ strongly fast invertible, we prove that there is no Lyapunov transformation bringing the system into diagonal form. For our particular fundamental matrix, it holds
		\begin{equation*}
			\left(\frac{\langle x_1,x_2\rangle}{\|x_1\|\,\|x_2\|}\right)^2=\frac{1}{1+\left(\frac{f}{g}\right)^{-2}}.
		\end{equation*}
		Since $f(t_{n,3})/g(t_{n,3})=1/\epsilon_n\to\infty$, the above expression converges to $1$ along the subsequence $(t_{n,3})_{n}$. Now, let $Y=[y_1,y_2]$ be an arbitrary fundamental matrix. By expanding  $y_1$ and $y_2$ with respect to our particular solutions and by using 
		$$\lim_{t\to\infty}\frac{\|x_2(t)\|}{\|x_1(t)\|}= 0,$$
		one computes
		\begin{align*}
			\lim_{n\to\infty}\left(\frac{\langle y_1(t_{n,3}),y_2(t_{n,3})\rangle}{\|y_1(t_{n,3})\|\,\|y_2(t_{n,3})\|}\right)^2= 1.
		\end{align*}
		Thus, 
		\begin{align*}
			\inf_t\frac{G(Y)}{G(y_1)G(y_2)}=1-\sup_t\left(\frac{\langle y_1,y_2\rangle}{\|y_1\|\,\|y_2\|}\right)^2=0,
		\end{align*}
		and \cref{thm:BlockTriang} implies that there is no Lyapunov transformation bringing the system into diagonal form. In particular, the system cannot be $1$-dim.\ strongly fast invertible according to \cref{thm:FICharacterization}.		
	\end{proof}

\subsubsection{Perturbation theory}

	We now establish several perturbation results for strongly fast invertible systems. Since we plan to use them to derive convergence theorems for Benettin's algorithm and since numerical integration can be represented as a piecewise continuous dynamical system, the perturbation results are formulated for systems
	\begin{equation}\label{eqn:PCSystem}
		\dot{x}=A(t)x
	\end{equation}
	such that $A$ is bounded and piecewise continuous. Moreover, we assume all perturbations in this subsection to be of the same class, i.e., bounded and piecewise continuous.\par 

	While strong fast invertibility is weaker than the stability of Lyapunov exponents, it is enough to ensure upper semicontinuity of $\lambda_{1,L}$.\par 

	\begin{theorem}\label[theorem]{thm:UpperSemicont}
		Assume \cref{eqn:PCSystem} is $L$-dim.\ strongly fast invertible. For any perturbation $Q(t)$, the largest Lyapunov exponent of the induced system satisfies
		\begin{align*}
			\tilde{\lambda}_{1,L}&\leq \lambda_{1,L} + Lc_{FI,L}\limsup_{t\to\infty}\frac{1}{t}\int_{0}^{t}\|Q(\tau)\|\,d\tau\\
			&\leq \lambda_{1,L} + Lc_{FI,L}\|Q\|_{L^{\infty}}
		\end{align*}
	\end{theorem}
	
	\begin{remark}
		Remember that
		$$\lambda_{1,L}=\Lambda_1+\dots+\Lambda_L$$
		for regular systems.
	\end{remark}

	\begin{proof}
		Let $\tilde{X}(t)$ be a fundamental matrix of the perturbed system. $\wedge^L\tilde{X}$ solves the induced equation
		\begin{align*}
			\frac{d}{dt}(\wedge^L\tilde{X})&=\widehat{A+Q}^L(\wedge^L\tilde{X})\\
			&=\hat{A}^L(\wedge^LX)+\hat{Q}^L(\wedge^LX)
		\end{align*}
		and can be seen as a solution of a perturbed system to
		$$\frac{d}{dt}(\wedge^LX)=\widehat{A}^L(\wedge^LX).$$
		The variation of constants formula implies
		\begin{align*}
			&\wedge^L\tilde{X}(t,s)= \wedge^LX(t,s)\\
			&\hspace{2em}+ \int_{s}^{t}(\wedge^LX(t,\tau))\hat{Q}^L(\tau)(\wedge^L\tilde{X}(\tau,s))\,d\tau.
		\end{align*}
		Using $\|\hat{Q}^L(t)\|\leq L\|Q(t)\|$, we first estimate 
		\begin{align*}
			&\frac{\|\wedge^L\tilde{X}(t,s)\|}{\|\wedge^LX(t,s)\|}\\
			&\hspace{1em}\leq 1 + L\int_{s}^{t}\frac{\|\wedge^LX(t,\tau)\|\,\|\wedge^LX(\tau,s)\|}{\|\wedge^LX(t,s)\|}\\
			&\hspace{2em}\cdot\|Q(\tau)\|\frac{\|\wedge^L\tilde{X}(\tau,s)\|}{\|\wedge^LX(\tau,s)\|}\,d\tau\\
			&\hspace{1em}\leq 1 + Lc_{FI,L}\int_{s}^{t}\|Q(\tau)\|\frac{\|\wedge^L\tilde{X}(\tau,s)\|}{\|\wedge^LX(\tau,s)\|}\,d\tau
		\end{align*}
		and then apply Grönwall's inequality to get
		\begin{equation*}
			\frac{\|\wedge^L\tilde{X}(t,s)\|}{\|\wedge^LX(t,s)\|}\leq e^{Lc_{FI,L}\int_{s}^{t}\|Q(\tau)\|\,d\tau}.
		\end{equation*}
		The claim now follows from \cref{prop:CharExpMatrixNorm}.
	\end{proof} 	
	
	\begin{lemma}
		Perturbations on finite intervals do not change the Lyapunov spectrum of a system or its induced systems.
	\end{lemma}
	
	\begin{proof}
		Since characteristic exponents are asymptotic quantities and since one may choose fundamental matrices $X(t)$ and $\tilde{X}(t)$ that coincide for large $t$, the claim follows.
	\end{proof}
	
	\begin{theorem}\label{thm:LowerEqual}
		Assume \cref{eqn:PCSystem} is $L$-dim.\ strongly fast invertible. If $Q(t)$ is a perturbation such that
		\begin{equation*}
			\|Q(t)\|\to 0\text{ as }t\to\infty,
		\end{equation*}
		then
		\begin{equation*}
			\tilde{\lambda}_{1,L}\leq \lambda_{1,L}.
		\end{equation*}
	\end{theorem}
	
	\begin{proof}
		The theorem is a direct consequence of \cref{thm:UpperSemicont}.\par 
	\end{proof}	
	
	Equality in \cref{thm:LowerEqual} can be achieved for $L^1$-perturbations.\par 

	\begin{theorem}\label{thm:Equal}
		Assume \cref{eqn:PCSystem} is $L$-dim.\ strongly fast invertible. If $Q(t)$ is a perturbation such that
		\begin{equation*}
			\int_0^{\infty}\|Q(t)\|\,dt<\infty,
		\end{equation*}
		then
		\begin{equation*}
			\tilde{\lambda}_{1,L} = \lambda_{1,L}.
		\end{equation*}
	\end{theorem}

	The theorem follows directly from \cref{thm:LowerEqual} using the following lemma.

	\begin{lemma}\label[lemma]{lemma:FICont}
		Assume \cref{eqn:PCSystem} is $L$-dim.\ strongly fast invertible. If $Q(t)$ is a perturbation such that
		\begin{equation*}
			\int_0^{\infty}\|Q(t)\|\,dt<\infty,
		\end{equation*}
		then the perturbed system is $L$-dim.\ strongly fast invertible. Moreover, $c_{FI,L}$ is continuous with respect to $L^1$-perturbations. 
	\end{lemma}

	\begin{proof}
		According to \cref{lemma:FIExtPowers} strong fast invertibility can be tested on $[T,\infty)$ for any fixed $T>0$. Thus, by setting $Q(t)$ to zero on a finite interval, we may assume that 
		\begin{equation}\label{eqn:IntBound}
			\int_0^{\infty}\|Q(t)\|\,dt<\frac{\log(2)}{Lc_{FI,L}}
		\end{equation}		
		without affecting whether the perturbed system is strongly fast invertible or not. The variation of constants formula implies
		\begin{align*}
			&\wedge^L\tilde{X}(t,s)-\wedge^LX(t,s)\\
			&\hspace{1em}=\int_{s}^{t}(\wedge^LX(t,\tau))\hat{Q}^L(\tau)(\wedge^L\tilde{X}(\tau,s))\,d\tau.
		\end{align*}
		Thus, it holds
		\begin{align*}
			&1+\frac{\|\wedge^L\tilde{X}(t,s)-\wedge^LX(t,s)\|}{\|\wedge^LX(t,s)\|}\\
			&\hspace{1em}\leq 1 + L\int_{s}^{t}\frac{\|\wedge^LX(t,\tau)\|\,\|\wedge^LX(\tau,s)\|}{\|\wedge^LX(t,s)\|}\|Q(\tau)\|\\
			&\hspace{2em}\cdot\left(1+\frac{\|\wedge^L\tilde{X}(\tau,s)-\wedge^LX(\tau,s)\|}{\|(\wedge^LX(\tau,s)\|}\right)\,d\tau.
		\end{align*}
		Through Grönwall's inequality we get
		\begin{align*}
			&1+\frac{\|\wedge^L\tilde{X}(t,s)-\wedge^LX(t,s)\|}{\|\wedge^LX(t,s)\|}\\
			&\hspace{1em}\leq e^{Lc_{FI,L}\int_{s}^{t}\|Q(\tau)\|\,d\tau}
		\end{align*}
		and due to \cref{eqn:IntBound}
		\small
		\begin{align*}
			&\|\wedge^L\tilde{X}(t,s)\|\\
			&\hspace{1em}\geq \|\wedge^LX(t,s)\|\left(1-\frac{\|\wedge^L\tilde{X}(t,s)-\wedge^LX(t,s)\|}{\|\wedge^LX(t,s)\|}\right)\\
			&\hspace{1em}\geq \|\wedge^LX(t,s)\|\left(2-e^{Lc_{FI,L}\int_{s}^{t}\|Q(\tau)\|\,d\tau}\right).
		\end{align*}
		\normalsize
		Combining the last estimate with the one from the proof of \cref{thm:UpperSemicont} yields
		\begin{align*}
			&\frac{\|\wedge^L\tilde{X}(t,\tau)\|\,\|\wedge^L\tilde{X}(\tau,s)\|}{\|\wedge^L\tilde{X}(t,s)\|}\\
			&\hspace{1em}\leq \frac{\|\wedge^LX(t,\tau)\|\,\|\wedge^LX(\tau,s)\|}{\|\wedge^LX(t,s)\|}\\
			&\hspace{2em}\cdot\frac{e^{Lc_{FI,L}\int_{s}^{t}\|Q(u)\|\,du}}{2-e^{Lc_{FI,L}\int_{s}^{t}\|Q(u)\|\,du}},
		\end{align*}
		which proves that the perturbed system is $L$-dim.\ strongly fast invertible.\par 
		
		Moreover, it follows that
		\begin{equation*}
			\tilde{c}_{FI,L}\leq c_{FI,L}\frac{e^{Lc_{FI,L}\|Q\|_{L^1}}}{2-e^{Lc_{FI,L}\|Q\|_{L^1}}},
		\end{equation*}
		which implies upper semicontinuity of $c_{FI,L}$ with respect to $L^1$-perturbations. Lower semi-continuity follows by switching roles, i.e., by viewing the original system as a perturbation of the perturbed system via $-Q(t)$.\par 
	\end{proof}
	
	\begin{corollary}\label[corollary]{cor:EqualRegular}
		Assume \cref{eqn:PCSystem} is regular and $L$-dim.\ strongly fast invertible. For every $\hat{v}\notin V_{2,L}'$, we find $\epsilon>0$ such that 
		$$\lim_{t\to\infty}\frac{1}{t}\log\|(\wedge^L\tilde{X}(t,0))\hat{v}\|=\Lambda_1+\dots+\Lambda_L$$		
		for any perturbation $Q(t)$ with
		\begin{equation*}
			\int_0^{\infty}\|Q(t)\|\,dt<\epsilon.
		\end{equation*}
	\end{corollary}
	
	\begin{proof}
		Since characteristic exponents are invariant under multiplication of the solution by a nonzero constant, we may assume that $\|\hat{v}\|=1$. It holds
		\small
		\begin{align*}
			&\|(\wedge^L \tilde{X}(t,0))\hat{v}\|\\
			&\hspace{1em}\geq \|(\wedge^LX(t,0))\hat{v}\|-\|(\wedge^L\tilde{X}(t,0)-\wedge^LX(t,0))\hat{v}\|\\
			&\hspace{1em}\geq \|(\wedge^LX(t,0))\hat{v}\|\\
			&\hspace{1.2em}\cdot\bigg(1-\frac{\|\wedge^L\tilde{X}(t,0)-\wedge^LX(t,0)\|}{\|\wedge^LX(t,0)\|}\frac{\|\wedge^LX(t,0)\|}{\|(\wedge^LX(t,0))\hat{v}\|}\bigg).
		\end{align*}
		\normalsize
		Since the system is regular, \cref{cor:StrongestGrowthInducedSystem} implies that
		$$\lim_{t\to\infty}\frac{\|\wedge^LX(t,0)\|}{\|(\wedge^LX(t,0))\hat{v}\|}<\infty.$$
		Now, according to the proof of \cref{lemma:FICont}, we may choose $\epsilon>0$ small enough to ensure 
		$$\|(\wedge^L \tilde{X}(t,0))\hat{v}\|\geq c\|(\wedge^LX(t,0))\hat{v}\|$$
		for some constant $c>0$ independent of $t$. Hence, it holds
		\begin{align*}
			&\liminf_{t\to\infty}\frac{1}{t}\log\|(\wedge^L\tilde{X}(t,0))\hat{v}\|\\
			&\hspace{1em}\geq \liminf_{t\to\infty}\frac{1}{t}\log\|(\wedge^LX(t,0))\hat{v}\|\\
			&\hspace{1em} = \lambda_{1,L}
		\end{align*}
		and 
		\begin{align*}
			&\limsup_{t\to\infty}\frac{1}{t}\log\|(\wedge^L\tilde{X}(t,0))\hat{v}\|\leq  \tilde{\lambda}_{1,L} =  \lambda_{1,L}
		\end{align*}
		due to \cref{thm:Equal}. The claim follows since regularity implies $ \lambda_{1,L}  = \Lambda_1+\dots+\Lambda_L$.
	\end{proof}

\bigskip

	\section{Computation of Lyapunov exponents}

	In this section we derive convergence results for the computation of Lyapunov exponents via \emph{Benettin's algorithm} \cite{BenettinEtAl1980LyapunovCharacteristicExponents, BenettinEtAl1980LyapunovCharacteristicExponentsa}.\par 
	
	Assume $h_{\max}<\infty$ and let $$\Phi:\mathbb{R}_{\geq 0}\times[0,h_{\max}]\to GL(d,\mathbb{R})$$ be a (linear) one-step method that is consistent of order $p>0$.
	
	\begin{definition}
		We call $\Phi$ \emph{consistent} if there is a constant $c_{\Phi}>0$ such that
		\begin{equation*}
			\|\Phi(t,h)-X(t+h,t)\|\leq c_{\Phi}h^{p+1}
		\end{equation*}
		for all $t\geq 0$ and $0\leq h\leq h_{\max}$.
	\end{definition}
	
	Given stepsizes $0<h_n\leq h_{\max}$, we shorten our notation by defining
	\begin{equation*}
		\Phi_n:=\Phi(t_{n-1}, h_n)
	\end{equation*}
	and
	\begin{equation*}
		\Phi^n :=\Phi_{n}\dots\Phi_{1}.
	\end{equation*}
	A similar notation will be adopted for $X$.\par

\subsection{Benettin's algorithm}

	The idea behind Benettin's algorithm can be explained via exterior products: Since $V_{2,L}$ is a proper subspace of $\wedge^L\mathbb{R}^d$, the solution to Lebesgue-almost every initial condition $v_1\wedge\dots\wedge v_L$ has characteristic exponent $\lambda_{1,L}$. Using $X(t)[v_1,\dots,v_L]=Q(t)R(t)$, it holds 
	\begin{align*}
		\chi[r_{11}(t)\dots r_{LL}(t)] &= \chi[(\wedge^LX(t))(v_1\wedge\dots\wedge v_L)] \\
		&= \lambda_{1,L}
	\end{align*}
	and hence, for regular systems, the Lyapunov exponents can be computed as
	\begin{equation*}
		\Lambda_i=\lim_{t\to\infty}\frac{1}{t}\log r_{ii}(t)
	\end{equation*}
	for Lebesgue-almost every tuple of initial vectors, which is the core idea behind Benettin's algorithm.\par 
	
	\begin{algorithm}
		\DontPrintSemicolon
		
		\KwInput{number of integration steps $N$, stepsizes $(h_n)_{n=1}^N$}
		\KwOutput{computed Lyapunov exponents $\mu_1,\dots,\mu_d$}
		
		$V_0 =$ rand(d)\tcp*{set random initial vectors}
		\For{n = 1:N}    
		{ 
			$W = \Phi_nV_{n-1}$\tcp*{evolve}
			$[Q_n,R_n] = \textnormal{qr}(W)$\tcp*{orthonormalize}
			$V_n = Q_n$\tcp*{set new vectors}
		}
		\For{i = 1:d}    
		{ 
			$\mu_i(N) = \frac{1}{t_N}\sum_{n=1}^{N}\log(R_n)_{ii}$\tcp*{average}
		}
		
		\caption{Benettin's algorithm \cite{BenettinEtAl1980LyapunovCharacteristicExponents, BenettinEtAl1980LyapunovCharacteristicExponentsa}}
	\end{algorithm}
	
	The propagated vectors\footnote{The propagated vectors from Benettin's algorithm are more than a mere byproduct for Lyapunov exponents. For instance, they have been exploited by Ginelli et al.\,\cite{GinelliEtAl2007CharacterizingDynamicsCovariant} and by Wolfe-Samelson \cite{WolfeSamelson2007EfficientMethodRecovering} in their algorithms to compute covariant Lyapunov vectors (see \cite{Noethen2019ProjectorbasedConvergenceProof,Noethen2021ComputingCovariantLyapunov} for a theoretical analysis).} (columns of $V_n$) are reorthonormalized periodically to prevent numerical singularities. If not, they could collapse onto the fastest expanding direction rendering them numerically indistinguishable, which makes it impossible to compute their associated volumes. Analytically, however, intermediate orthonormalizations do not play a role since the associated volumes stay the same. Indeed, we get the same output analytically if we perform the $QR$-decomposition only once at the end:
	\begin{align*}
		&\Phi^nV_0=\Phi_n\dots\Phi_2\Phi_1V_0=\Phi_n\dots\Phi_2V_1R_1\\
		&\hspace{1em}=\dots = V_n (R_n\dots R_1).
	\end{align*}
	
	One may adjust the frequency of orthonormalizations depending on how fast $V_n$ becomes singular. Furthermore, in practice it makes sense to compute the Lyapunov exponents as a running average during the propagation loop to save memory and to monitor convergence properties.\par 
	
	While the algorithm works perfectly fine analytically, there are some numerical challenges. In practice several types of errors influence the output of Benettin's algorithm. For example, errors are introduced by numerical integration, by limiting the integration time, by finite-precision computing or through errors of the underlying equations\footnote{If the linear system is derived as the linearization along a solution of a nonlinear system, integration errors from computing the solution of the nonlinear system can result in errors in the derived linear equations of order $O(1)$. In that case our perturbation results for linear theory are ineffective and nonlinear theory is necessary to study the computational error.}.\par

\subsection{Convergence results}	
	
	For all results in this subsection, we assume:
	\begin{itemize}
		\item $0< h_n\leq h_{\max}$,
		\item $t_n=h_1+\dots+h_n\to\infty$,
		\item $\Phi$ is consistent of order $p>0$.
	\end{itemize}
	
	Before deriving the desired convergence theorems for Benettin's algorithm, we introduce two auxiliary systems: a piecewise constant approximation of our original system and a piecewise constant system representing the numerical integration.\par 

	If $\|X_n-I\|<1$, the logarithm of $X_n$ exists and is equal to
	\begin{equation*}
		\log X_n = \sum_{k=1}^{\infty}(-1)^{k+1}\frac{(X_n-I)^k}{k}.
	\end{equation*}	
	We set
	\begin{equation*}
		A_n:=\frac{1}{h_n}\log X_n.
	\end{equation*}
	Since 
	\begin{align*}
		\|X_n-I\|&\leq \int_{t_{n-1}}^{t_n}\|A(s)X(s,t_{n-1})\|\,ds\\
		&\leq Me^{h_{\max}M}h_n,
	\end{align*}
	\normalsize
	$A_n$ is well-defined for small stepsizes.\par 
	
	Similarly, we set
	\begin{equation*}
		B_n:=\frac{1}{h_n}\log\Phi_n
	\end{equation*}
	whenever $\|\Phi_n-I\|<1$. Consistency of the numerical integrator and the previous estimate imply that $B_n$ is well-defined for small stepsizes as well.\par 
	
	We have the following relation between $A_n$ and $B_n$:\par 
	 
	\begin{lemma}\label[lemma]{lemma:EstimatePiecewiseSystems}
		There is a constant\,\footnote{A more detailed analysis of the constant for systems with stable Lyapunov exponents can be found in \cite{DiecivanVleck2005ErrorComputingLyapunov}.} $c>0$ independent of the stepsizes such that
		$$\|A_n-B_n\|\leq c h_n^p$$
		for small $h_n$.
	\end{lemma}
	
	\begin{proof}
		Assume $h_n$ is small enough such that
		\begin{equation}\label{eqn:HBarRequirement}
			\max(\|\Phi_n-I\|,\|X_n-\Phi_n\|)\leq \frac{1}{4}.
		\end{equation}
		T he logarithms of $X_n$ and $\Phi_n$ exist. Moreover, it holds
		\begin{align*}
			&\|\log(X_n)-\log(\Phi_n)\|\\
			&\hspace{1em}\leq \sum_{k=1}^{\infty}\frac{\|(X_n-I)^k-(\Phi_{n}-I)^k\|}{k}.
		\end{align*}
		Using 
		$$\|(M+E)^k-M^k)\|\leq \|E\|2^k\max(\|M\|,\|E\|)^{k-1}$$
		with $M:=\Phi_n-I$ and $E:=X_n-\Phi_n$, we estimate
		\begin{align*}
			&\|A_n-B_n\|\\
			&\hspace{1em}=\frac{1}{h_n}\|\log(X_n)-\log(\Phi_n)\|\\
			&\hspace{1em}\leq \frac{1}{h_n}\left(\sum_{k=1}^{\infty}\frac{4}{2^kk}\right)\|X_n-\Phi_n\|.
		\end{align*}
		The claim follows from consistency of the numerical integrator.
	\end{proof}
	
	We may regard $A_n$ as a perturbation of the original system on $[t_{n-1},t_n)$. The size of the perturbation depends on the continuity property of $A$.\par 
	
	\begin{lemma}
		If $A$ is globally Lipschitz continuous, then there is a constant $c>0$ independent of the stepsizes such that
		$$\sup_{t\in[t_{n-1},t_n)}\|A(t)-A_n\|\leq ch_n$$
		for small $h_n$.
	\end{lemma}
	
	\begin{proof}
		If $h_n$ is small enough such that \cref{eqn:HBarRequirement} holds, then
		$$\|X_n-I\|\leq \frac{1}{2}$$
		and
		\begin{align*}
			&\|\log(X_n)-(X_n-I)\|\\
			&\hspace{1em}\leq\bigg(\sum_{k=2}^{\infty}\frac{\|X_n-I\|^{k-2}}{k}\bigg)\|X_n-I\|^2\\
			&\hspace{1em}\leq \bigg(\underset{=:C}{\underbrace{\sum_{k=2}^{\infty}\frac{4}{2^{k}k}\bigg)M^2e^{2h_{\max}M}}}h_n^2.
		\end{align*}
		Thus, it holds
		\begin{align*}
			&\|A(t)-A_n\|\\
			&\hspace{1em}=\left\|A(t)-\frac{1}{h_n}\log X_n\right\|\\
			&\hspace{1em}\leq\left\|A(t)-\frac{1}{h_n}(X_n-I)\right\|+Ch_n\\
			&\hspace{1em}\leq \frac{1}{h_n}\int_{t_{n-1}}^{t_n}\|A(t)-A(s)X(s,t_{n-1})\|\,ds+ Ch_n.
		\end{align*}
		for $t\in[t_{n-1},t_n)$. Now, the claim follows from
		\begin{align*}
			&A(t)-A(s)X(s,t_{n-1})=A(t)-A(s)\\
			&\hspace{2em}+ A(s)(X(t_{n-1},t_{n-1})-X(s,t_{n-1}))
		\end{align*}	
		via Lipschitz continuity of $A$.
	\end{proof}
	
	In the following, we only consider stepsizes such that $h_n\to 0$. While we cannot guarantee that $A_n$ and $B_n$ are well-defined for small $n$, there is $N\geq 0$ such that $A_n$ and $B_n$ are well-defined for $n\geq N$. We define
	$$A_{pc}(t):=\begin{cases}
		0,\quad& t\in [0,t_{N-1})\\
		A_n,\quad& t\in[t_{n-1},t_n)\text{ for } n\geq N
	\end{cases}$$
	and 
	$$B_{pc}(t):=\begin{cases}
		0,\quad&t\in [0,t_{N-1})\\
		B_n.\quad& t\in[t_{n-1},t_n)\text{ for } n\geq N
	\end{cases}$$
	By possibly increasing $N$, we ensure that $A_{pc}$ and $B_{pc}$ are bounded and can be estimated using the previous lemmata.\footnote{This is the case if \cref{eqn:HBarRequirement} is satisfied for $n\geq N$.}\par 
	
	Let us denote the Cauchy matrix corresponding to $A_{pc}$ by $X_{pc}(t,s)$. If $n\geq m\geq N-1$, then
	\begin{equation*}
		X_{pc}(t_n,t_m)=e^{h_nA_n}\dots e^{h_{m+1}A_{m+1}}=X(t_n,t_m).
	\end{equation*}
	Hence, $A_{pc}$ can be seen as a piecewise constant approximation of our original system (for large $t$).\par  
	
	\begin{lemma}
		The Lyapunov spectra of the original system and the piecewise constant approximation coincide. The same holds for their induced systems. Moreover, if the original system is regular, then the piecewise constant approximation is regular.
	\end{lemma}
	
	\begin{proof}
		Since 
		\begin{align*}
			{X}_{pc}(t,t_{N-1})= e^{(t-t_{n-1})A_n} X(t,t_{n-1})^{-1} X(t,t_{N-1})
		\end{align*}
		for $t\in[t_{n-1},t_n)$ and $n\geq N$, the growth rates of solutions with the same initial condition differ by at most a constant:		
		\begin{align*}
			\left(\frac{\|{X}_{pc}(t,t_{N-1})v\|}{\|X(t,t_{N-1})v\|}\right)^{\pm 1}\leq e^{h_{\max}(\|{A}_{pc}\|_{\infty}+M)}
		\end{align*}
		for any $v\neq 0$. Hence, the Lyapunov spectra of the original system and the piecewise constant approximation coincide. The remaining statements follow similarly.
	\end{proof}
	
	We write $\tilde{X}_{pc}(t,s)$ for the Cauchy matrix corresponding to $B_{pc}$ and note that
	$$\tilde{X}_{pc}(t_n,t_m)=\Phi_n\dots\Phi_{m+1}$$
	for $n\geq m\geq N-1$. In particular, fixing the fundamental matrix $\tilde{X}_{pc}(t)$ satisfying $\tilde{X}_{pc}(t_{N-1})=\Phi^{N-1}$, we have
	\begin{align*}
		\tilde{X}_{pc}(t_n)&=\tilde{X}_{pc}(t_n,t_{N-1})\tilde{X}_{pc}(t_{N-1})\\
		&=\Phi_n\dots\Phi_{N}\Phi^{N-1}\\
		&=\Phi^n
	\end{align*}
	for $n\geq N-1$. Thus, $B_{pc}$ can be understood as representing the numerical integration. It follows that
	\begin{align*}
		&\mu_1(n)+\dots+\mu_L(n)\\
		&\hspace{1em}=\sum_{i=1}^{L}\frac{1}{t_n}\sum_{k=1}^{n}\log(R_k)_{ii}\\
		&\hspace{1em}=\frac{1}{t_n}\log\left(\prod_{i=1}^{L}(R_n\dots R_1)_{ii}\right)\\
		&\hspace{1em}=\frac{1}{t_n}\log \|(\wedge^L\Phi_n)(v_1\wedge\dots\wedge v_L)\|\\
		&\hspace{1em}=\frac{1}{t_n}\log \|(\wedge^L\tilde{X}_{pc}(t_n))(v_1\wedge\dots\wedge v_L)\|
	\end{align*}
	for $n\geq N-1$. Since the stepsizes and ${B}_{pc}$ are bounded, we have
	\begin{align*}
		&\limsup_{n\to\infty}\mu_1(n)+\dots+\mu_L(n)\\
		&\hspace{1em}=\limsup_{n\to\infty}\frac{1}{t_n}\log \|(\wedge^L\tilde{X}_{pc}(t_n))(v_1\wedge\dots\wedge v_L)\|\\
		&\hspace{1em}=\limsup_{t\to\infty}\frac{1}{t}\log \|(\wedge^L\tilde{X}_{pc}(t))(v_1\wedge\dots\wedge v_L)\|.
	\end{align*}
	Thus, Benettin's algorithm computes $\tilde{\lambda}_{1,L}$ of ${B}_{pc}$ for Lebesgue-almost every tuple of initial vectors:
	\begin{equation}\label{eqn:BenettinOutputLimSup}
		\limsup_{n\to\infty}\mu_1(n)+\dots+\mu_L(n)=\tilde{\lambda}_{1,L}.
	\end{equation}
	In particular, if ${B}_{pc}$ is regular, then
	\begin{equation}\label{eqn:BenettinOutputLim}
		\lim_{n\to\infty}\mu_1(n)+\dots+\mu_L(n)=\tilde{\Lambda}_1+\dots+\tilde{\Lambda}_L.
	\end{equation}
	
	\begin{theorem}\label[theorem]{thm:BenettinStableLEsLowEqual}
		Assume \cref{eqn:system} is globally Lipschitz continuous and has stable Lyapunov exponents. If $h_n\to 0$, then 
		$$\limsup_{n\to\infty}\mu_1(n)+\dots+\mu_i(n)\leq \Lambda_1+\dots+\Lambda_i$$
		for all $i$ and Lebesgue-almost every tuple of initial vectors.\footnote{A version of the estimate for finite time can be found in \cite{DiecivanVleck2006PerturbationTheoryApproximation}. The authors relate the error to the departure from normality of the $R$-matrices obtained during Benettin's algorithm.}
	\end{theorem}
	
	\begin{proof}
		Since there are constant $c_1,c_2>0$ independent of the stepsizes such that
		\begin{align*}
			&\sup_{t\in[t_{n-1},t_n)}\|A(t)-{B}_{pc}(t)\|\\
			&\hspace{1em}\leq \sup_{t\in[t_{n-1},t_n)}\|A(t)-A_n\|+\|A_n-B_n\|\\
			&\hspace{1em}\leq c_1h_n+c_2h_n^p
		\end{align*}
		for $n\geq N$, ${B}_{pc}$ is a perturbation of the original system such that $\|Q(t)\|\to0$. In particular, \cref{thm:PerturbationToZero} implies that ${B}_{pc}$ has the same Lyapunov spectrum as the original system. Now, \cref{eqn:BenettinOutputLimSup} and \cref{prop:LE1Induced} imply the theorem.
	\end{proof}
	
	\begin{theorem}\label[theorem]{thm:BenettinStableLEsEqual}
		Assume \cref{eqn:system} is globally Lipschitz continuous, regular and has stable Lyapunov exponents. If $h_n\to 0$, then
		$$\lim_{n\to\infty}\mu_i(n)=\Lambda_i$$
		for all $i$ and Lebesgue-almost every tuple of initial vectors.
	\end{theorem}
	
	\begin{proof}
		The claim follows similarly to the last theorem using \cref{prop:PerturbationToZeroRegular} and \cref{eqn:BenettinOutputLim} instead of \cref{eqn:BenettinOutputLimSup}.
	\end{proof}
	
	\begin{remark}
		If the stability of Lyapunov exponents transfers from $A$ to ${A}_{pc}$, then the assumption that $A$ is globally Lipschitz continuous can be dropped in \cref{thm:BenettinStableLEsLowEqual,thm:BenettinStableLEsEqual}.
	\end{remark}
	
	Next, we derive convergence results for strongly fast invertible systems. For this, we need the following lemma:\par 
	 
	 \begin{lemma}
	 	The original system is $L$-dim.\ strongly fast invertible if and only if its piecewise constant approximation is $L$-dim.\ strongly fast invertible.
	 \end{lemma}
	 
	 \begin{proof}
	 	The claim follows from
	 	\begin{align*}
	 		&e^{(t_n-t)A_n}{X}_{pc}(t,s)e^{(s-t_{m-1})A_m}\\
	 		&\hspace{1em}=X_{pc}(t_{n},t_{m-1})\\
	 		&\hspace{1em}=X(t_{n},t_{m-1})\\
	 		&\hspace{1em}=X(t_n,t)X(t,s)X(s,t_{m-1})
	 	\end{align*}
	 	for $s\in[t_{m-1},t_m)$ and $t\in[t_{n-1},t_n)$, $n\geq  m\geq N$, and the fact that strong fast invertibility can be tested on $[t_{N-1},\infty)$ (see \cref{lemma:FIExtPowers}).
	 \end{proof}
	 
	 \begin{theorem}\label[theorem]{thm:BenettinFILEsEstimate}
	 	Assume \cref{eqn:system} is $L$-dim.\ strongly fast invertible. If $h_n\to 0$, then 
	 	$$\limsup_{n\to\infty}\mu_1(n)+\dots+\mu_L(n)\leq \Lambda_1+\dots+\Lambda_L$$
	 	for Lebesgue-almost every tuple of initial vectors.
	 \end{theorem}
	 
	 \begin{proof}
	 	Since ${A}_{pc}$ is $L$-dim.\ strongly fast invertible and its induced systems have the same Lyapunov spectra as the induced systems of the original system, \cref{thm:LowerEqual} and \cref{prop:LE1Induced} imply 
	 	$$\tilde{\lambda}_{1,L}\leq \lambda_{1,L}\leq \Lambda_1+\dots+\Lambda_L.$$
	 	The claim follows from \cref{eqn:BenettinOutputLimSup}.
	 \end{proof}
	 
	 \begin{theorem}\label[theorem]{thm:BenettinFILEsEqual}
	 	Assume \cref{eqn:system} is regular and $L$-dim.\ strongly fast invertible. If
	 	$$\sum_{n=1}^{\infty} h_n^{p+1}<\infty,$$
	 	then
	 	\begin{equation*}
	 		\limsup_{n\to\infty}\mu_1(n)+\dots+\mu_L(n)= \Lambda_1+\dots+\Lambda_L
	 	\end{equation*}
	 	for Lebesgue-almost every tuple of initial vectors.
	 \end{theorem}
	 
	 \begin{proof}
	 	\cref{lemma:EstimatePiecewiseSystems} and the stepsize condition ensure that 
	 	\begin{align*}
	 		&\int_{0}^{\infty}\|{A}_{pc}(t)-{B}_{pc}(t)\|\,dt<\infty.
	 	\end{align*}
	 	Now, the proof is as in \cref{thm:BenettinFILEsEstimate} using \cref{thm:Equal} and regularity:
	 	$$\tilde{\lambda}_{1,L}= \lambda_{1,L}= \Lambda_1+\dots+\Lambda_L.$$
	 \end{proof}
	 
	 Since we do not know if regularity transfers form the original to the numerical system, \cref{thm:BenettinFILEsEqual} only ensures convergence to the Lyapunov exponents as a limes superior. However, when fixing a tuple of initial vectors, we can ensure convergence as a limit if the stepsizes decay fast enough.\par  
	 
	 \begin{theorem}\label[theorem]{thm:BenettinFILEsEqualLim}
		 Assume \cref{eqn:system} is regular and $L$-dim.\ strongly fast invertible. For Lebesgue-almost every tuple of initial vectors, we have the following: If\ \ $\sum_{n=1}^{\infty} h_n^{p+1}$ is small enough, then
		 \begin{equation*}
		 	\lim_{n\to\infty}\mu_1(n)+\dots+\mu_L(n)= \Lambda_1+\dots+\Lambda_L.
		 \end{equation*}
	 \end{theorem}
	 
	 \begin{proof}
	 	The theorem follows from \cref{cor:EqualRegular}.
	 \end{proof}
	 
	 \begin{remark}
	 	The condition ``$\sum_{n=1}^{\infty} h_n^{p+1}$ small enough'' in \cref{thm:BenettinFILEsEqualLim} depends on the chosen initial vectors. Indeed, following the associated proofs, we require smaller stepsizes the smaller the first principle angle between $\textnormal{span}(v_1,\dots, v_L)$ and $V_{l+1}'$. 
	 \end{remark}
	 
	 Finally, when applied to the right dimensions, \cref{thm:BenettinFILEsEqual,thm:BenettinFILEsEqualLim} ensure that we may approximate the Lyapunov exponents of strongly fast invertible systems using Benettin's algorithm. 
	 
	 \begin{corollary}
	 	Assume \cref{eqn:system} is regular and strongly fast invertible at dim.\ $d_1+\dots+d_l$ for $l=1,\dots,p$.  If
	 	$$\sum_{n=1}^{\infty} h_n^{p+1}<\infty,$$
	 	then
	 	\begin{align*}
	 		d_i\lambda_i &= \limsup_{n\to\infty}\mu_1(n)+\dots+\mu_{d_1+\dots+d_i}(n)\\
	 		&\hspace{2em}-\limsup_{n\to\infty}\mu_1(n)+\dots+\mu_{d_1+\dots+d_{i-1}}(n)
	 	\end{align*}
	 	for all $i$ and for Lebesgue-almost every tuple of initial vectors. In particular, if the Lyapunov spectrum is simple, then
	 	\begin{align*}
	 		\lambda_i &= \limsup_{n\to\infty}\mu_1(n)+\dots+\mu_{i}(n)\\
	 		&\hspace{2em}-\limsup_{n\to\infty}\mu_1(n)+\dots+\mu_{i-1}(n)
	 	\end{align*}
	 	for all $i$ and for Lebesgue-almost every tuple of initial vectors.
	 \end{corollary}
	 
	 \begin{corollary}
		 Assume \cref{eqn:system} is regular and strongly fast invertible at dim.\ $d_1+\dots+d_l$ for $l=1,\dots,p$. For Lebesgue-almost every tuple of initial vectors, we have the following: If\ \ $\sum_{n=1}^{\infty} h_n^{p+1}$ is small enough, then
		 $$d_i\lambda_i = \lim_{n\to\infty}\mu_{d_1+\dots+d_{i-1}+1}(n)+\dots+ \mu_{d_1+\dots+d_i}(n)$$
		 for all $i$. In particular, if the Lyapunov spectrum is simple, then
		 $$\lambda_i=\lim_{n\to\infty}\mu_i(n)$$
		 for all $i$.
	 \end{corollary}

\bigskip

	\newpage

\section*{Acknowledgments}	
	
	This paper is a contribution to the project M1 (Dynamical Systems Methods and Reduced Models in Geophysical Fluid Dynamics) of the Collaborative Research Centre TRR 181 ``Energy Transfers in Atmosphere and Ocean'' funded by the Deutsche Forschungsgemeinschaft (DFG, German Research Foundation) - Projektnummer 274762653.
	
%
	
\bigskip

	\pagestyle{empty}


\printbibliography

\end{document}